\documentclass[11pt]{article}
\usepackage[a4paper]{geometry}  
\usepackage{times} 
\usepackage{xurl} 
\usepackage[italian,english]{babel}
\usepackage{latexsym} 
\usepackage{xcolor}

\makeatletter
\newcommand{\@BIBLABEL}{\@emptybiblabel}
\newcommand{\@emptybiblabel}[1]{}
\makeatother
\usepackage{graphicx}
\usepackage{latexsym}
\usepackage{amssymb}
\usepackage{bm}
\usepackage{amsmath}
\usepackage{epsfig}
\usepackage{caption}
\usepackage{subcaption}
\usepackage{graphicx}
\usepackage{caption}
\usepackage{subcaption}
\usepackage{mwe} 
\usepackage{morefloats}
\usepackage{cite}
\usepackage{epstopdf}
\usepackage[colorlinks=true,linkcolor=blue,citecolor=blue,urlcolor=blue]{hyperref}
\title{A Novel Gaussian filter-based Pressure Correction Technique with Super Compact Scheme for Unsteady 3D Incompressible, Viscous Flows}

\author{\textbf{Ashwani Punia$^{1}$, Rajendra K. Ray$^{2}$} \\
  1,2. School of Mathematical and Statistical Sciences, Indian Institute of Technology Mandi,\\
  Mandi, Himachal Pradesh, 175005, India \\
 {\tt mr.punia11@gmail.com}, {\tt rajendra.iitmandi.ac.in}}

\date{}

\begin{document}

\maketitle
\begin{abstract}
This work deals with a novel Gaussian filter-based pressure correction technique with a super compact higher order finite difference scheme for solving unsteady three-dimensional (3D) incompressible, viscous flows. This pressure correction technique offers significant advantages in terms of optimizing computational time by taking minimum iterations to reach the required accuracy, making it highly efficient and cost-effective. Pressure fields often exhibit highly nonlinear behavior, and employing the Gaussian filter can help to enhance their reliability by reducing noise and uncertainties. On the other hand, the super compact scheme uses minimum grid points to produce second-order accuracy in time and fourth-order accuracy in space variables. The main focus of this study is to enhance the accuracy and efficiency and minimize the computational cost of solving complex fluid flow problems. The super compact scheme utilizes 19 grid points at the known time level (i.e., $n^{th}$ time level) and only seven grid points from the unknown time level (i.e., $n + 1$ time level). By employing the above strategies, it becomes possible to notably decrease computational expenses while maintaining the accuracy of the computational scheme for solving complex fluid flow problems. We have implemented our methodology across three distinct scenarios: the 3D Burger's equation having analytical solution, and two variations of the lid-driven cavity problem. The outcomes of our numerical simulations exhibit a remarkable concordance, aligning exceptionally well with both the analytical benchmarks and previously validated numerical findings for the cavity problems.
This robust agreement underscores the efficacy and reliability of our approach. Furthermore, we conducted a detailed analysis of the pressure correction technique in terms of streamline contours and pressure correction results. The findings of this study have significant implications for various engineering and scientific disciplines, such as aerodynamics, hydrodynamics, and fluid-structure interaction analysis.

\end{abstract}
\vspace{0.4cm}
\par
\vspace{8pt}

\section{Introduction}
Numerical simulations have a crucial role in understanding complex fluid flow phenomena and predicting their behavior in various engineering applications. The Navier-Stokes equations, governing the dynamics of fluid motion, are fundamental in such simulations. These developments have occurred across three prominent methodologies: finite element, finite volume, and finite difference methods. While finite difference methods excel in their ease of implementation, both finite element and finite volume approaches offer greater flexibility when it comes to directly applying them to irregular domains. Nonetheless, it's worth acknowledging that every approach presents a unique array of benefits and hurdles.
Notably, the finite difference approach stands out as the most straightforward to implement. Over the past few years, a multitude of specialized methods have been pioneered, e.g., high-order compact (HOC) finite difference (FD) schemes known for their computational efficiency. Numerous researchers have innovated with fourth-order compact finite difference schemes, specifically designed for convection-diffusion equations on both two-dimensional uniform grids \cite{Karaa_2002, Zhang_2003, Wang_2010,Ray_2010} and three-dimensional spaces \cite{Ge_2002, Gupta_2000, Zhang_1998, Kalita_2014}. The finite difference approach relies primarily on Taylor series expansion, a powerful mathematical tool used to approximate functions. This simplicity in implementation has contributed to its popularity among researchers and engineers seeking to tackle complex fluid flow problems. By discretizing the governing equations into a grid-like structure, finite difference methods enable the conversion of continuous partial differential equations into discrete algebraic equations that can be numerically solved on a computer. It is important to emphasize that most finite difference schemes developed to address the Navier-Stokes equations have been tailored primarily for the context of two-dimensional (2D) flows. Notably, these schemes, extensively applied in the stream function-vorticity ($\psi-\omega$) formulation, have been the focus of studies conducted by researchers such as Ghia et al. (1982)\cite{Ghia_1982}, Lecointe and Piquet (1984)\cite{Lecointe_1984}, Weinan and Liu (1996) \cite{Weinan_1996}, Kupferman et al. (2001) \cite{Kupferman_2001}, Tian et al. (2003)\cite{Tian_2003}, Kalita and Chhabra (2006)\cite{Kalita_2006},  Ray et al. (2010) \cite{Ray_2010}, Wang et al. (2019)\cite{Wang_2019} and Yadav et al. (2023)\cite{Yadav_2023}. Consequently, due to their inherent limitations, these finite difference schemes were unable to be readily extended to handle three-dimensional (3D) flows. Only a small number of cases exist where high-order compact (HOC) schemes were specifically designed for three-dimensional convection-diffusion equations. However, even in these instances, the application scope remained restricted to scenarios involving steady-state linear convection-diffusion or the solution of 3D Poisson equations (Gupta and Kouatchou (1998)\cite{Gupta_1998}, Zhang et al. (2000)\cite{Zhang_2000a,Zhang_2000b}, McTaggart et al. (2004)\cite{McTaggart_2004}). The challenges in extending these HOC schemes to more complex 3D flows have remained a subject of ongoing research and exploration in the field of numerical simulations.\\
However, solving these equations can be computationally intensive, especially for three - dimensional problems or at high Reynolds numbers, which require fine grids and numerous time steps. To address these challenges, researchers have explored various numerical methods and techniques to enhance the efficiency and accuracy of simulations. One such approach is the Modified Compressibility technique\cite{Cortes_1994}, wherein pressure is iteratively solved to fulfill the incompressibility constraint within the Navier-Stokes equations. Although widely used, this technique can still impose significant computational costs, especially for large-scale problems. Kalita \cite{Kalita_2014} introduced an exceptionally efficient super compact higher-order scheme for simulating 3D incompressible flow. This approach utilized a modified compressibility method for pressure calculation, which, while effective, did have a significant drawback in terms of computational cost when applied to solving 3D fluid flow problems. It requires lots of pressure iterations at each time level, resulting in heightened computational costs.\\ 
In today's rapidly evolving computational simulation landscape, optimizing efficiency without compromising accuracy is a paramount goal. Reducing computational costs is crucial, especially in the context of simulating 3D fluid flow problems, where the computational demands are significantly high. In this context, we propose a new and efficient pressure correction technique that reduces the computational cost while maintaining high accuracy. 
Our approach is inspired by the need to accelerate the convergence of pressure iterations in the Modified Compressibility technique and enhance the overall efficiency of numerical simulations. By introducing a Gaussian-based pressure correction technique to the pressure field, we aim to achieve faster convergence rates and lower computational costs without compromising accuracy. The Gaussian filter is a commonly used image processing technique employed to enhance or modify images by reducing noise and smoothing out details while preserving the overall structure. This observation serves as inspiration for applying Gaussian filtering techniques in the context of pressure fields, where its capacity to reduce errors and fluctuations is particularly valuable. 
In this study, we present our proposed pressure correction technique and its application to three different fluid flow problems. Firstly, we investigate the validation of the three-dimensional Burger's equation, a fundamental nonlinear partial differential equation that serves as an ideal test case due to its known analytical solutions. We conduct a comparison between our computed results and the analytical solutions to assess the accuracy and efficiency of the scheme. Next, we apply our novel technique to two variations of the lid-driven cavity problem, a classic benchmark in fluid dynamics. The 2D lid-driven square cavity has been extensively studied, but its 3D counterpart remains relatively unexplored due to its increased complexity and memory requirements. Our approach allows us to efficiently tackle this challenging 3D problem at different Reynolds numbers.
Throughout our investigation, we systematically analyze the computational efficiency and importance of our proposed technique. We conduct a comparison of outcomes obtained through varying time steps and Reynolds numbers to demonstrate the robustness and versatility of our approach. Additionally, we investigate the flow patterns and pressure contours in the lid-driven cavity to gain deeper insights into the fluid behavior under different conditions.
In summary, this study aims to contribute to the advancement of numerical simulations in fluid dynamics by proposing a novel and efficient pressure correction technique. The success of our approach in accurately capturing flow phenomena while reducing computational costs has significant implications for a wide range of applications, from fundamental research to engineering simulations.
\section{Mathematical Modeling and Discretization}
\begin{figure}
    \centering
      \includegraphics[width=0.8\textwidth]{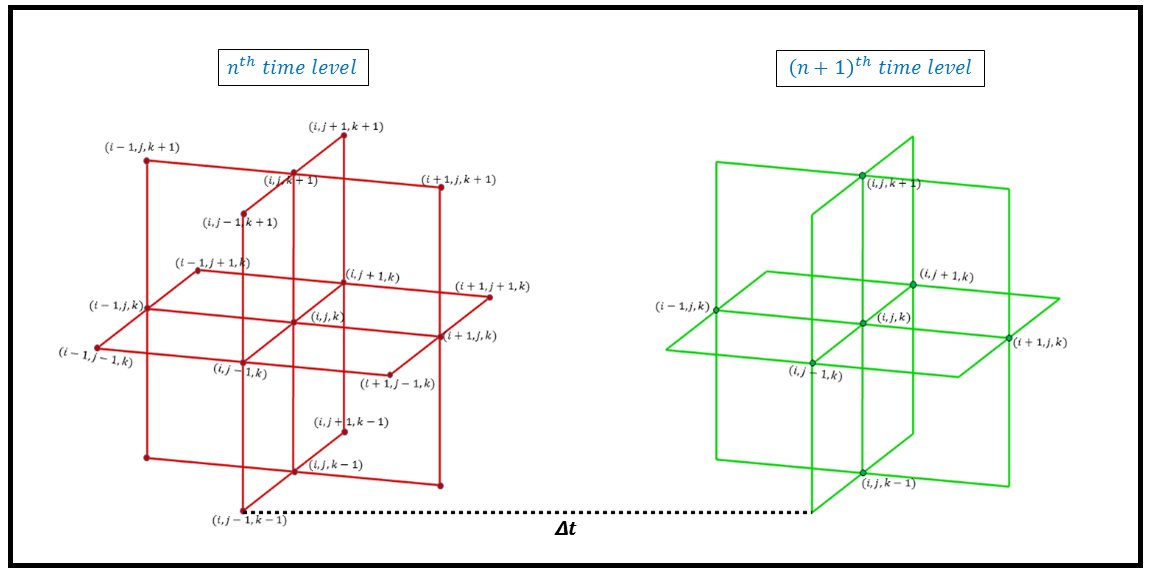}
    \caption{The super-compact unsteady stencil }
    \label{fig:stencil}
\end{figure}

The equation that governs the behavior of an unsteady three-dimensional (3D) convection-diffusion-reaction process for a transport variable denoted as ``$U$"  within a continuous domain, characterized by variable coefficients, can be expressed as follows.
\begin{equation}
\begin{aligned}
& m \frac{\partial U}{\partial t}+n(x, y, z, t) \frac{\partial U}{\partial x}+o(x, y, z, t) \frac{\partial U}{\partial y}+p(x, y, z, t) \frac{\partial U}{\partial z}+q(x, y, z, t) U \\
& \quad=\nabla^2 U+r(x, y, z, t),
\end{aligned}\label{eq1}
\end{equation}
In this context, the equation introduces a constant ``$m$" and several coefficients, namely ``$n$", ``$o$", and ``$p$" representing convection, ``$q$" denoting reaction, and ``$r$" indicating a forcing function. This equation models the convection-diffusion process of various fluid variables, including mass, heat, energy, vorticity, etc., within the given continuous domain. By appropriately selecting the values of ``$m, n, o, p, q$", and ``$r$", the equation can effectively represent the complete Navier-Stokes ($\mathrm{N}-\mathrm{S}$) equations. This equation thus stands as a unifying framework capable of describing an array of fluid dynamics phenomena within a singular mathematical construct.
It is essential to establish suitable boundary conditions for the domain to ensure a well-defined and physically meaningful problem formulation. Suppose the problem domain is cubical, and to discretize it, the mesh is characterized by incremental measures denoted as $h$, $k$, and $l$ along the $x$, $y$, and $z$ directions, respectively. This meticulous mesh creation allows us to systematically represent the domain in a structured manner. The super compact higher order finite difference scheme was recently developed by Kalita\cite{Kalita_2014}. Within this discretized framework, the Forward-Time Centered-Space (FTCS) scheme is applied to approximate the equation (\ref{eq1}) at the node $(i, j, k)$ of the mesh. The FTCS scheme is a numerical method commonly used to approximate time-dependent partial differential equations in computational fluid dynamics and other scientific simulations. It involves updating the solution at each node based on its current value and the values of neighboring nodes in a central differencing manner, where both time and space derivatives are approximated.
The conventional FTCS approximation for Eq. (\ref{eq1}) at the node $(i, j, k)$ can be expressed as follows:
\begin{equation}
\left(h \delta_t^{+}+f+p\delta_z+o \delta_y+n \delta_x-\delta_z^2-\delta_y^2-\delta_x^2\right) U_{i j k}-\xi_{i j k}=r_{i j k},
\end{equation}\label{eq2}
In the given equation, the function $U$ is defined on a three-dimensional grid with grid points $(x_i, y_j, z_k)$ denoted as $U_{i j k}$. The mathematical apparatus introduced includes operators that play crucial roles in this context. These operators encompass $\delta_x, \delta_x^2, \delta_y, \delta_y^2, \delta_z, \delta_z^2,$ and $\delta_t^{+}$. Each of these operators corresponds to specific differentiations: first- and second-order central differences along the $x$-, $y$-, and $z$-directions, and a first-order forward difference along the temporal direction, respectively.  The truncation error $\xi_{i j k}$ for this numerical approach, with a uniform time step $\Delta t$, quantifies the error introduced due to this discretization is given by
\begin{equation}
\begin{aligned}
\xi_{i j k}= & {\left[-\frac{h^2}{12}\left(\frac{\partial^4 U}{\partial x^4}-2 n \frac{\partial^3 U}{\partial x^3}\right)-\frac{k^2}{12}\left(\frac{\partial^4 U}{\partial y^4} - 2 o \frac{\partial^3 U}{\partial y^3}\right)\right.} \\
& \left.-\frac{l^2}{12}\left(\frac{\partial^4 U}{\partial z^4}-2 p \frac{\partial^3 U}{\partial z^3}\right)+m \frac{\Delta t}{2} \frac{\partial^2 U}{\partial t^2}\right]_{i j k}+O\left(\Delta t^2, h^4, k^4, l^4\right) .
\end{aligned}\label{eq3}
\end{equation}
In order to attain a higher level of precision in time (second-order accuracy) as well as spatial precision (fourth-order accuracy) for the equation (\ref{eq1}), the derivatives of the leading term in equation (\ref{eq3}) are approximated in a compact manner \cite{MacKinnon_1991, Spotz_1995}, leading to a formulation with reduced truncation error. To achieve the objective, the initial partial differential equation (PDE) given by Eq. (\ref{eq1}) is handled as an additional relation from which higher derivatives can be derived. To illustrate, the forward temporal difference method for the transport variable $U$ and the backward difference method for variables $m, n, o, p, q $, $r$\cite{Kalita_2014} are utilized. This allows us to express the derivative in the first term on the right-hand side of Eq. (\ref{eq3}) as follows: 
\begin{equation}
\begin{aligned}
\left.m \frac{\partial^2 U}{\partial t^2}\right|_{i j k}= & \left(\delta_z^2+\delta_y^2+\delta_x^2-q_{i j k}-n_{i j k} \delta_x-o_{i j k} \delta_y-p_{i j k} \delta_z\right) \delta_t^{+} U_{i j k} \\
& -\left(\delta_t^{-} n_{i j k} \delta_x+\delta_t^{-} o_{i j k} \delta_y+\delta_t^{-} q_{i j k}+\delta_t^{-} p_{i j k} \delta_z +\delta_t^{-} r_{i j k}\right) U_{i j k}+O\left(\Delta t, h^2, k^2, l^2\right),
\end{aligned} \label{eq4}
\end{equation}
The operator $\delta_t^{-}$ represents a first-order backward difference in the temporal domain. Likewise, we can devise analogous approximations for the spatial derivatives as well.
Consequently, by substituting the derivatives in Eq. (\ref{eq3}) with the approximations given in Eq. (\ref{eq4}) and similar expressions, and then replacing $\xi_{i j k}$ in Eq. (\ref{eq2}) accordingly, we obtain the following approximation\cite{Kalita_2014} of order $O\left(\Delta t^2, h^4, k^4, l^4\right)$ for our main governing equation (\ref{eq1}). 
$$
\begin{aligned}
m[1 & +\left(\frac{h^2}{12}-\frac{\Delta t}{2 m}\right)\left(\delta_x^2-n_{i j k} \delta_x\right)+\left(\frac{k^2}{12}-\frac{\Delta t}{2 m}\right)\left(\delta_y^2-o_{i j k} \delta_y\right) \\
& \left.+\left(\frac{l^2}{12}-\frac{\Delta t}{2 m}\right)\left(\delta_z^2-p_{i j k} \delta_z\right)+\frac{\Delta t}{2 m} q_{i j k}\right] \delta_t^{+} U_{i j k} \\
& +\left(-\alpha_{i j k} \delta_x^2-\beta_{i j k} \delta_y^2-\gamma_{i j k} \delta_z^2+A_{i j k} \delta_x+B_{i j k} \delta_y+C_{i j k} \delta_z+D_{i j k}\right) U_{i j k} \\
& -\frac{h^2+k^2}{12}\left(\delta_x^2 \delta_y^2-n_{i j} \delta_x \delta_y^2-o_{i j k} \delta_x^2 \delta_y-p1_{i j k} \delta_x \delta_y\right) U_{i j k}
\end{aligned}
$$

\begin{equation}
\begin{aligned}
& -\frac{k^2+l^2}{12}\left(\delta_y^2 \delta_z^2-o_{i j k} \delta_y \delta_z^2-p_{i j k} \delta_y^2 \delta_z-q1_{i j k} \delta_y \delta_z\right) U_{i j k} \\
& -\frac{l^2+h^2}{12}\left(\delta_z^2 \delta_x^2-p_{i j k} \delta_z \delta_x^2-n_{i j k} \delta_z^2 \delta_x-r1_{i j k} \delta_z \delta_x\right) U_{i j k} \\
= &  R_{i j k} .
\end{aligned}\label{eq5}
\end{equation}
The coefficients $\alpha_{i j k}, \beta_{i j k}, \gamma_{i j k}, A_{i j k}, B_{i j k}, C_{i j k}, D_{i j k}, R_{i j k}, p1_{i j k}, q1_{i j k}$ and $r1_{i j k}$ are as follows:
$$
\begin{aligned}
& \alpha_{i j k}=\frac{h^2}{12}\left(n_{i j k}^2-q_{i j k}-2 \delta_x n_{i j k}\right)+1 \text {, } \\
& \beta_{i j k}=\frac{k^2}{12}\left(o_{i j k}^2-q_{i j k}-2 \delta_y o_{i j k}\right)+1 \text {, } \\
& \gamma_{i j k}=\frac{l^2}{12}\left(p_{i j k}^2-q_{i j k}-2 \delta_z p_{i j k}\right)+1\text {, } \\
& A_{i j k}=\left[\frac{h^2}{12}\left(\delta_x^2-n_{i j k} \delta_x\right)+\frac{k^2}{12}\left(\delta_y^2-o_{i j k} \delta_y\right)+\frac{l^2}{12}\left(\delta_z^2-p_{i j k} \delta_z\right)+\frac{\Delta t}{2} \delta_t^{-}+1\right] n_{i j k} \\
& -\frac{h^2}{12}\left(n_{i j k}-2 \delta_x\right) q_{i j k}, \\
& B_{i j k}=\left[\frac{h^2}{12}\left(\delta_x^2-n_{i j k} \delta_x\right)+\frac{k^2}{12}\left(\delta_y^2-o_{i j k} \delta_y\right)+\frac{l^2}{12}\left(\delta_z^2-p_{i j k} \delta_z\right)+\frac{\Delta t}{2} \delta_t^{-}+1\right] o_{i j k} \\
& -\frac{k^2}{12}\left(o_{i j k}-2 \delta_y\right) q_{i j k} \text {, } \\
& C_{i j k}=\left[\frac{h^2}{12}\left(\delta_x^2-n_{i j k} \delta_x\right)+\frac{k^2}{12}\left(\delta_y^2-o_{i j k} \delta_y\right)+\frac{l^2}{12}\left(\delta_z^2-p_{i j k} \delta_z\right)+\frac{\Delta t}{2} \delta_t^{-}+1\right] p_{i j k} \\
& -\frac{l^2}{12}\left(p_{i j k}-2 \delta_z\right) q_{i j k} \text {, } \\
& D_{i j k}=\left[\frac{h^2}{12}\left(\delta_x^2-n_{i j k} \delta_x\right)+\frac{k^2}{12}\left(\delta_y^2-o_{i j k} \delta_y\right)+\frac{l^2}{12}\left(\delta_z^2-p_{i j k} \delta_z\right)+\frac{\Delta t}{2} \delta_t^{-}+1\right] q_{i j k}, \\
& R_{i j k}=\left[\frac{h^2}{12}\left(\delta_x^2-n_{i j k} \delta_x\right)+\frac{k^2}{12}\left(\delta_y^2-o_{i j k} \delta_y\right)+\frac{l^2}{12}\left(\delta_z^2-p_{i j k} \delta_z\right)+\frac{\Delta t}{2} \delta_t^{-}+1\right] r_{i j k}, \\
& p1_{i j k}=-n_{i j k} o_{i j k}+\frac{2}{h^2+k^2}\left(k^2 \delta_y n_{i j k}+h^2 \delta_x o_{i j k}\right) \text {, } \\
& q1_{i j k}=-o_{i j k} p_{i j k}+\frac{2}{k^2+l^2}\left(l^2 \delta_z o_{i j k}+k^2 \delta_y p_{i j k}\right) \text {, } \\
& r1_{i j k}=-p_{i j k} n_{i j k}+\frac{2}{l^2+h^2}\left(h^2 \delta_x p_{i j k}+l^2 \delta_z n_{i j k}\right) . \\
&
\end{aligned}
$$
With this (\ref{eq5}), an implicit finite difference scheme is obtained that exhibits second-order accuracy in time \(O(t^2)\), fourth-order accuracy \(O(h^4, k^4, l^4)\) in spatial grid spacing and achieved using a $(19,7)$ stencil, as depicted in Figure \ref{fig:stencil}.

In this formulation (as shown in Eq. (\ref{eq4})), a good approach is employed to avoid mixed derivative terms \cite{Kalita_2014}. This approach results in a compact seven-point stencil at the (n + 1)$^{th}$ time step, significantly reducing the computational burden. Notably, many high-order compact schemes developed for two-dimensional \cite{Kalita_2002,Ray_2017} convection-diffusion equations typically necessitate nine points in the stencil at the (n + 1)$^{th}$ time level. The super compact higher order compact scheme offers double advantages. Firstly, it allows for a simpler seven-point stencil, requiring only the $(i, j, k)$ $^{th}$ point and its six neighboring points (as shown in Figure \ref{fig:stencil}) at the (n + 1)$^{th}$ time level. Secondly, it effectively eliminates the need for extensive corner points and considerably reduces the number of points needed for the approximation, enhancing computational efficiency. The more information regarding the higher-order super compact scheme can be found in the reference \cite{Kalita_2014}.

\subsection{Utilizing the Navier-Stokes ($\mathrm{N}-\mathrm{S}$) Equations}
We can express the transient three-dimensional incompressible Navier-Stokes ($\mathrm{N}-\mathrm{S}$) equations in a non-dimensional form as follows:
\begin{equation}
\begin{aligned}
\frac{\partial U}{\partial t}+U \cdot \nabla U & =-\nabla pr+\frac{1}{R e} \nabla^2 U,
\end{aligned}\label{eqn_1_1} 
\end{equation}
\begin{equation}
   \begin{aligned}
\nabla \cdot U & =0,
\end{aligned} \label{eqn_1_2} 
\end{equation}
The given equations describe the transient three-dimensional incompressible Navier-Stokes (N-S) equations in non-dimensional form. In this representation, the velocity vector is denoted by $U=(u, v, w)$, where $u$, $v$, and $w$ represent the velocity components in the $x$, $y$, and $z$ directions, respectively. The variable $t$ corresponds to time, and $pr$ represents the pressure. To non-dimensionalize the equations, the Reynolds number ($Re$) is introduced, which is defined as $Re=\frac{U_0 N}{v^*}$, where $N$ is a characteristic length, $U_0$ is a characteristic velocity, and $v^*$ is the kinematic viscosity.
The operator $\nabla$ represents the gradient operator and is given by $\nabla=\frac{\partial}{\partial x} \mathbf{i}+\frac{\partial}{\partial y} \mathbf{j}+\frac{\partial}{\partial z} \mathbf{k}$. Equation (\ref{eqn_1_1}) pertains to the momentum equation, while equation (\ref{eqn_1_2}) represents the continuity equation. These equations can also be expressed in terms of the primitive variables as follow.
\begin{equation}
\begin{aligned}
\frac{\partial u}{\partial t}-\frac{1}{R e}\left(\frac{\partial^2 u}{\partial z^2}+\frac{\partial^2 u}{\partial y^2}+\frac{\partial^2 u}{\partial x^2}\right)=-v \frac{\partial u}{\partial y}-u \frac{\partial u}{\partial x}-w \frac{\partial u}{\partial z}-\frac{\partial pr}{\partial x}, 
\end{aligned}\label{eqn_1_3} 
\end{equation}

\begin{equation}
\frac{\partial v}{\partial t}-\frac{1}{R e}\left(\frac{\partial^2 v}{\partial z^2}+\frac{\partial^2 v}{\partial y^2}+\frac{\partial^2 v}{\partial x^2}\right)=-v \frac{\partial v}{\partial y}-u \frac{\partial v}{\partial x}-w \frac{\partial v}{\partial z}-\frac{\partial pr}{\partial y}, \label{eqn_1_4} 
\end{equation}

\begin{equation}
\frac{\partial w}{\partial t}-\frac{1}{R e}\left(\frac{\partial^2 w}{\partial z^2}+\frac{\partial^2 w}{\partial y^2}+\frac{\partial^2 w}{\partial x^2}\right)=-v \frac{\partial w}{\partial y}-u \frac{\partial w}{\partial x}-w \frac{\partial w}{\partial z}-\frac{\partial pr}{\partial z},\label{eqn_1_5}  
\end{equation}

\begin{equation}
\frac{\partial v}{\partial y}+\frac{\partial u}{\partial x}+\frac{\partial w}{\partial z}=0.\label{eqn_1_6}
\end{equation} 
To solve the Navier-Stokes equations, the first step involves discretizing the momentum equations (Eqs. (\ref{eqn_1_3}),(\ref{eqn_1_4}),(\ref{eqn_1_5})) using a specific numerical scheme denoted as Eq. (\ref{eq5}). In this discretization process, we assign the values of parameters as follows:\\
$m$ is set to the Reynolds number ($Re$).\\
$n$ is set to $Re \cdot u$, where $u$ represents the velocity component in the $x$ direction.\\
$o$ is set to $Re \cdot v$, where $v$ represents the velocity component in the $y$ direction.\\
$p$ is set to $Re \cdot w$, where $w$ represents the velocity component in the $z$ direction.\\
$q$ is assigned a value of 0.\\
$r$ represents one of three options: $Re \cdot \frac{\partial pr}{\partial x}$, $Re \cdot \frac{\partial pr}{\partial y}$, or $Re \cdot \frac{\partial pr}{\partial z}$, depending on the particular derivative term in the momentum equations. By employing this discretization approach, the Navier-Stokes equations can be transformed into a system of algebraic equations that can be solved efficiently using some advanced iterative solver. To solve the system of algebraic equations, we use the Biconjugate Gradient Stabilized method.\\ 
It should be emphasized that analytical expressions for the pressure variable ($pr$) are not available. Therefore, we need to resort numerical approximations for the pressure gradients. In this regard, central difference approximations are often used for interior points, while standard one-sided first or higher-order approximations are employed at the boundaries.
Once the momentum equations (Eqs. (\ref{eqn_1_3}),(\ref{eqn_1_4}),(\ref{eqn_1_5})) are solved using the super compact higher order numerical finite difference scheme, the next step is to determine the pressure field ($pr$). However, this task becomes challenging as the incompressible Navier-Stokes equations lack an explicit pressure term. This complexity arises from the choice of the primitive variable representation, which necessitates considering the pressure variable during the computational process. So, the main point is that the solution of the incompressible Navier-Stokes equations demands careful handling of the pressure term, requiring numerical approximations and diligent consideration of the primitive variables throughout the calculations.\\
In the subsequent section, we address the challenge of dealing with the pressure variable in a highly effective and cost-efficient manner. We present a comprehensive approach that allows us to handle the pressure-related complexities while ensuring computational efficiency and accuracy. This method not only enables us to obtain reliable solutions for the incompressible Navier-Stokes equations but also significantly reduces the computational resources required for the simulations. By leveraging this innovative pressure correction technique, we can overcome the absence of an explicit pressure term and efficiently incorporate the pressure variable into our calculations, leading to a robust and practical solution strategy.
\section{Proposed Pressure Correction Technique}
This section introduces an innovative approach for managing pressure calculation within the solution of the Navier-Stokes equations, employing the modified compressibility technique. In general, handling the pressure equation in the Navier-Stokes equation is the most challenging task for any numerical method due to its highly non-linear behavior. To overcome this issue, the modified compressibility technique\cite{Cortes_1994} is an efficient and widely used approach for simulating fluid flows, where the pressure equation is solved iteratively to satisfy the incompressibility constraint. Though this technique is quite efficient, but it takes quite a large number of iterations to reach the desired accuracy, which increases the computational cost. Our proposed pressure correction technique aims to further enhance computational efficiency by reducing the number of iterations required to the desired accuracy of the pressure solution.
\subsection{Existing Pressure Iteration Approach}
In the method proposed by Cortes and Miller (1994) \cite{Cortes_1994}, they introduced a modification to the continuity equation to tackle the pressure-related challenges. The modified continuity equation takes the form:
$$
pr+\lambda \nabla \cdot \mathbf{v}=0 .
$$
During each time step, after calculating the pressure gradients and solving the momentum equations, they compute the dilation parameter $D$ which can be calculated as $u_x+v_y+w_z$. If the maximum absolute value of $D$, denoted as $|D|{\max}$, is below a predefined tolerance limit, the pressure value is supposed to reach the desired accuracy, and overall computation proceeds to the next time level, and follows the same steps again. Nonetheless, in the event that the maximum value of $|D|$ surpasses the specified tolerance limit, they initiate a pressure correction step is introduced to refine the pressure value:
$$
pr^{\text {new }}=pr^{\text {old }}-\lambda \nabla \cdot \mathbf{v} .
$$
where $p^{\text {new}}$ is the updated pressure, $p^{\text {old }}$ is the pressure value obtained in the previous iteration, and $\lambda$ is a relaxation parameter. The process is repeated iteratively until the maximum absolute value of the divergence of the velocity field, denoted as $|\nabla \cdot \mathbf{v}|_{\max}$, meets the tolerance limit.

\subsection{New Pressure Correction Technique}
The proposed pressure correction technique aims to accelerate the convergence of the pressure calculation process at each time level. Instead of relying solely on the conventional pressure iteration approach, we introduce an additional step after calculating the pressure at each pressure iteration. During each time step, after calculating the pressure gradients and solving the momentum equations using the super compact higher-order scheme, we apply a Gaussian smoothing function to the pressure field. The Gaussian smoothing function is a well-known mathematical tool that can effectively reduce noise and fluctuations in data. It is especially effective in image processing and computer vision applications to decrease noise and highlight key features. In image processing, the core concept behind a Gaussian filter is convolution. In this process, a small matrix, known as a kernel or filter, is applied to each pixel in the image. This filter gives more importance to nearby pixels while progressively reducing the impact of distant ones, creating a smoothing effect and removing the noise from the image. 
This observation serves as inspiration for applying Gaussian filtering techniques in the context of pressure fields, where its capacity to reduce errors and fluctuations is particularly valuable. Pressure fields often exhibit highly nonlinear behavior, and employing the Gaussian filter can help to enhance their reliability by reducing noise and uncertainties.\\ 
In our approach, we convolve the pressure field with a Gaussian kernel to reduce the error as well as any rapid variations that might hinder the convergence process. The main advantage of the proposed pressure correction technique lies in its ability to accelerate convergence and significantly reduce the number of iterations required to reach a solution. By applying the Gaussian smoothing function to the pressure field, we dampen high-frequency oscillations, enabling faster convergence without sacrificing accuracy. Reducing the number of pressure iterations translates to substantial reduction in computational time, especially for large-scale 3D simulations and complex flow problems. The proposed technique can significantly improve the efficiency of the modified compressibility technique, making it even more attractive for practical engineering applications. To validate the efficacy of our proposed technique, we conducted extensive numerical studies on various benchmark test cases. The results demonstrate a notable reduction in the number of iterations required to calculate the pressure value of the desired accuracy, maintaining accurate and stable solutions.\\
\begin{center}
\textit{\textbf{Algorithm:}} Proposed Pressure Handling Technique
\end{center}
\begin{enumerate}
\item Initialize the pressure field $pr^{\text{new}} = pr^{\text{old}}$.
\item Calculate the pressure gradients and solve the momentum equations using the Super Compact higher order scheme.
\item Increment the loop for pressure iterations.
\item Compute the divergence of the velocity field $|\nabla \cdot \mathbf{v}|_{\max}$.
\item If $|\nabla \cdot \mathbf{v}|_{\max}$ $>$ tolerance, do the following:\\
(a) Apply the pressure correction: $p^{\text{new}} = G_\xi * p^{\text{old}} - \lambda \nabla \cdot \mathbf{v}$, where $\lambda$ is a relaxation parameter and the former pressure field is convoluted with the Gaussian kernel to contribute to the updated solution.\\
(b) Go to step 2\\
Else, go to step 6.
\item The pressure field reaches the desired accuracy, go to the next time level.\\
\end{enumerate}
The proposed pressure correction technique utilizes the Modified Compressibility technique along with the Gaussian smoothing function to accelerate the convergence of the pressure iteration process. Additionally, It ensures that the continuity equation is effectively satisfied. By reducing noise (error) and rapid variations in the pressure field, the technique achieves faster convergence while maintaining the accuracy of the  solutions in each time step. Indeed, the implementation of our proposed pressure handling approach offers several advantages over traditional techniques such as the Pressure-Poisson Equation (PPE) approach and artificial compressibility methods. Notably, it achieves faster convergence and computational efficiency while maintaining simplicity and straightforwardness in its application.

\section{Results and Discussion}
\subsection{Problem-1}
We first consider the three-dimensional Burger’s equation, which has analytical solution. Fletcher (1983)\cite{Fletcher_1983} derived the steady solutions of the two-dimensional Burgers' equations and subsequently, $\mathrm{Xu}$ et al. (1997)\cite{Xu_1997} extended this approach to three dimensions in time. The governing equations are expressed as follows:
\begin{equation}
\frac{\partial U}{\partial t}-\frac{1}{R e}\left(\frac{\partial^2 U}{\partial x^2}+\frac{\partial^2 U}{\partial y^2}+\frac{\partial^2 U}{\partial z^2}\right)=-u \frac{\partial U}{\partial x}-v \frac{\partial U}{\partial y}-w \frac{\partial U}{\partial z},
\label{eq_2_1}
\end{equation}
Here, $U$ signifies the three-dimensional velocity field, while $u$, $v$, and $w$ refer to the velocity components along the $x, y,$ and $z$ directions, respectively. The parameter $Re$ corresponds to the Reynolds number, and it influences the behavior of the flow. This equation governs the evolution of fluid flow in three dimensions, and it plays a crucial role in understanding various fluid phenomena. For the purpose of code validation and to ensure the accuracy of the super compact scheme \cite{Kalita_2014}, we have chosen this problem. The solutions provided below are analytical solutions to the nonlinear equations (\ref{eq_2_1}):
$$
\begin{aligned}
u & =-\frac{2d_0}{R e}\left[c_2 e^{-\lambda\left(\frac{t}{R e}\right)} n_x \pi \cos \left(n_x \pi x\right) \sin \left(n_y \pi y\right) \sin \left(n_z \pi z\right)\right]  \\
v & =-\frac{2d_0}{R e}\left[c_2 e^{-\lambda\left(\frac{t}{R e}\right)} n_x \pi \cos \left(n_y \pi y\right) \sin \left(n_x \pi x\right) \sin \left(n_z \pi z\right)\right] \\
w & =-\frac{2d_0}{R e}\left[c_2 e^{-\lambda\left(\frac{t}{R e}\right)} n_x \pi \cos \left(n_z \pi z\right) \sin \left(n_x \pi x\right) \sin \left(n_y \pi y\right)\right] 
\end{aligned}
$$
where
$$
\begin{aligned}
\lambda & =\pi^2\left(n_x^2+n_y^2+n_z^2\right) \\
d_0 & =\frac{1}{c_1+c_2 e^{-\lambda\left(\frac{t}{R e}\right)} \sin \left(n_x \pi x\right) \sin \left(n_y \pi y\right) \sin \left(n_z \pi z\right)}
\end{aligned}
$$

\begin{figure}[htbp]
 \centering
 \vspace*{5pt}%
 \hspace*{\fill}%
\begin{subfigure}{0.50\textwidth}     
    \centering
    \includegraphics[width=\textwidth]{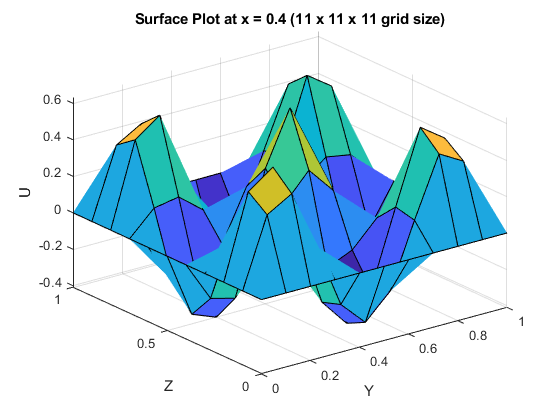}%
    \captionsetup{skip=5pt}%
    \caption{(a)}
    \label{fig:surface_plot_U_x0.5_11_Exact}
  \end{subfigure}%
 \begin{subfigure}{0.50\textwidth}        
   \centering
    \includegraphics[width=\textwidth]{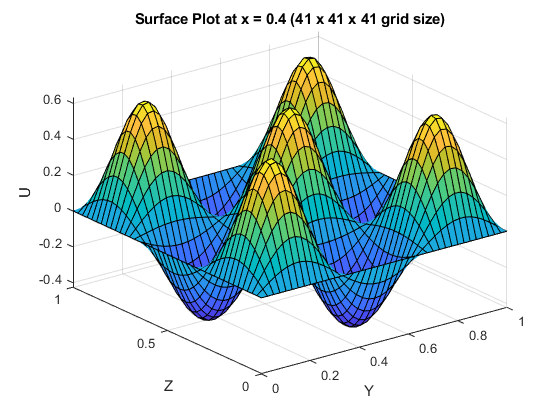}%
    \captionsetup{skip=5pt}%
    \caption{(b)}
    \label{fig:surface_plot_U_x0.5_61_NUM}
  \end{subfigure}
  \hspace*{\fill}

  \vspace*{8pt}%
  \hspace*{\fill}%

  \begin{subfigure}{0.50\textwidth}     
    \centering
    \includegraphics[width=\textwidth]{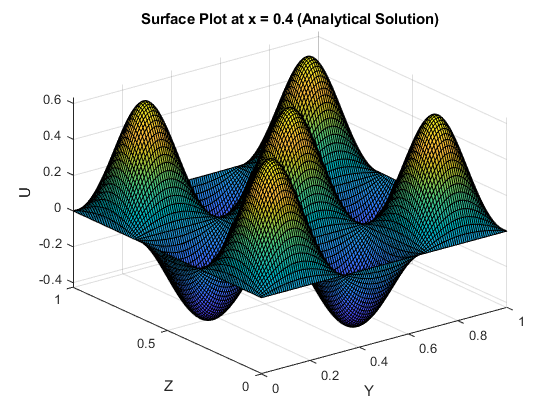}%
    \captionsetup{skip=5pt}%
    \caption{(c)}
    \label{fig:surface_plot_U_x0.5_61_ANALYTICAL}
  \end{subfigure}%
  \caption{Comparison of surface Plots: Numerical vs. Analytical Solution for $Re = 10$ at $t = 1.0$ and $x=0.4$ plane (a) Numerical solution using $11\times 11\times11$ grid size (b) Numerical solution using $41\times 41\times41$ grid size (c) Analytical solution}
  \label{fig:problem1}
\end{figure}


\begin{figure}[htbp]
 \centering
 \vspace*{5pt}%
 \hspace*{\fill}%
\begin{subfigure}{0.50\textwidth}     
    \centering
    \includegraphics[width=\textwidth]{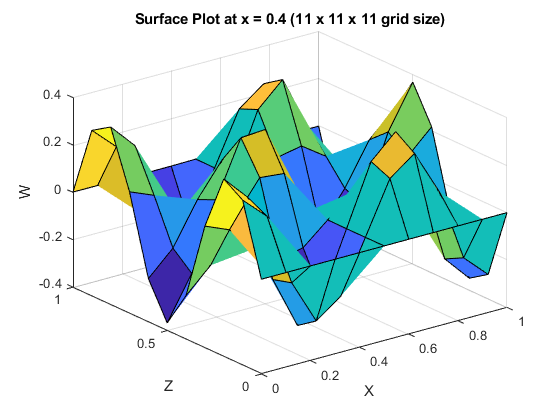}%
    \captionsetup{skip=5pt}%
    \caption{(a)}
    \label{fig:surface_plot_U_x0.5_11_Exact_W}
  \end{subfigure}%
 \begin{subfigure}{0.50\textwidth}        
   \centering
    \includegraphics[width=\textwidth]{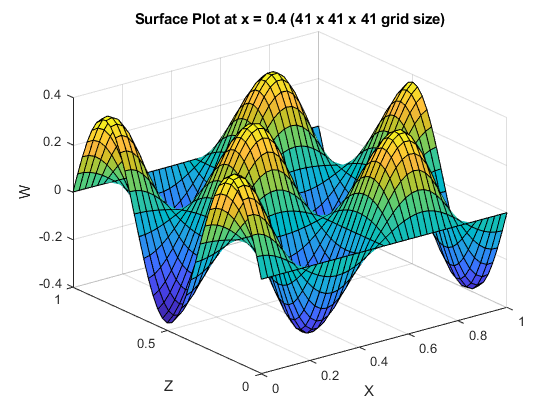}%
    \captionsetup{skip=5pt}%
    \caption{(b)}
    \label{fig:surface_plot_U_x0.5_61_NUM_W}
  \end{subfigure}
  \hspace*{\fill}

  \vspace*{8pt}%
  \hspace*{\fill}%
  
  \begin{subfigure}{0.50\textwidth}     
    \centering
    \includegraphics[width=\textwidth]{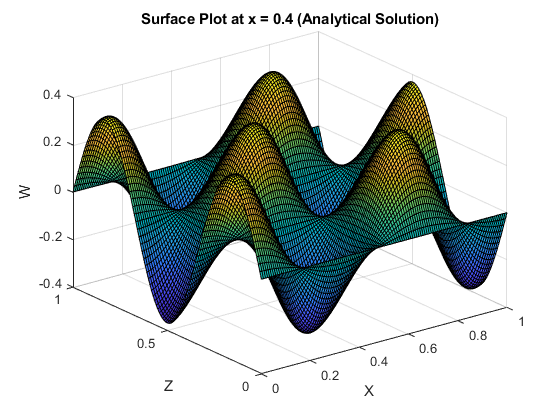}%
    \captionsetup{skip=5pt}%
    \caption{(c)}
    \label{fig:surface_plot_U_x0.5_61_Analyticall_W}
  \end{subfigure}%
  \caption{Comparison of Surface Plots: Numerical vs. Analytical Solution for $Re = 10$ at $t = 1.0$ and $x=0.4$ plane (a) Numerical solution using $11\times 11\times11$ grid size (b) Numerical solution using $41\times 41\times 41$ grid size (c) Analytical solution}
  \label{fig:problem1_W}

\end{figure}
The constants $c_1$ and $c_2$ are adjustable parameters that govern the amplitude of these analytical solutions.
The comparison between surface plots of the analytical and present computed solutions is presented in Figures \ref{fig:problem1} and \ref{fig:problem1_W} for two different grid sizes: ($11\times11\times11$) and ($41\times41\times41$), focusing on the $x=0.4$ and $y=0.4$ planes, respectively. For our comparisons, we have considered the following parameter values in the analytical solution. $3.0 = n_z =  n_y = n_x $, $c_1=6.0$, $c_2=2$ and $t=1.0$ for $Re=10$. From Figures \ref{fig:problem1} and \ref{fig:problem1_W}, one can observe that the difference between the analytical and numerical solutions is insignificant. The numerical outcomes are in excellent agreement with the analytical solutions, demonstrating the accuracy of the super compact scheme. To study the order of accuracy of the method, different computational errors are calculated and are presented in Table \ref{rate_convergence} for three different grid sizes: $11 \times 11 \times 11$, $21 \times 21 \times 21$, and $41 \times 41 \times 41$, at time, $t = 1.0$, $5.0$, with $Re = 10$ and time step, $\Delta t = 0.01$. The table provides the Average absolute errors, Average relative errors, and root mean square relative errors on each grid size.
The definition of the root mean square relative error\cite{Kalita_2014} is as follows:
$$
E=\frac{R M S\left(U_{\text {ana }}-U_{\text {num }}\right)}{R M S\left(U_{\text {ana }}\right)}
$$
Here, $R M S\left(U_{\text {ana }}-U_{\text {num }}\right)$ represents the root mean square of the differences between the values of $U_{\text {ana }}$ and $U_{\text {num}}$ taken at each node and $R M S\left(U_{\text {ana }}\right)$ is the root mean square value of the analytical solution. Upon careful analysis, it becomes evident that the method is fourth-order accurate, which confirms the theoretical result. More details about the super compact scheme and its accuracy can be seen in the reference \cite{Kalita_2014}.

{\small\begin{table}[htbp]
\caption{\small Comparison of Convergence Rates: Average Absolute errors ($e_1$), Average Relative errors ($e_2$), and Average Root Mean Square errors ($e_3$) at Different Grid Sizes and Time $ t = 1.0, 5.0 $}\label{rate_convergence}
\centering
 \begin{tabular}{cccc}  \hline \hline
& M$\times$N$\times$P    &   Error    &  Rate of Convergence    \\ \hline 
$t=1.0$\\
&(11 $\times$ 11 $\times$ 11)           &  9.0130E-6  &    -   \\
$e_1$ &(21 $\times$ 21 $\times$ 21)         &  6.1534E-7   & 3.873  \\
&(41 $\times$ 41 $\times$ 41)             &  4.2314E-8  &  3.860     \\
\hline
&(11 $\times$ 11 $\times$ 11)           &  3.471E-4  &    -   \\
$e_2$ &(21 $\times$ 21 $\times$ 21)         &  2.3293E-5   & 3.897  \\
&(41 $\times$ 41 $\times$ 41)             &  1.571E-6  &  3.890     \\
\hline
&(11 $\times$ 11 $\times$ 11)           &  7.223E-4  &    -   \\
$e_2$ &(21 $\times$ 21 $\times$ 21)         &  4.171E-5   & 4.114  \\
&(41 $\times$ 41 $\times$ 41)             &  2.610E-6  &  3.998     \\
\hline\hline
$t=5.0$\\
&(11 $\times$ 11 $\times$ 11)           &  5.8674E-5  &    -   \\
$e_1$ &(21 $\times$ 21 $\times$ 21)         &  3.8931E-6   & 3.921  \\
&(41 $\times$ 41 $\times$ 41)             &  2.5953E-7  &  3.900     \\
\hline
&(11 $\times$ 11 $\times$ 11)           &  4.7241E-3  &    -   \\
$e_2$ &(21 $\times$ 21 $\times$ 21)         &  3.1172E-4   & 3.921  \\
&(41 $\times$ 41 $\times$ 41)             &  2.0872E-5  &  3.900     \\
\hline
&(11 $\times$ 11 $\times$ 11)           &  1.3122E-2  &    -   \\
$e_2$ &(21 $\times$ 21 $\times$ 21)         &  7.9814E-4   & 4.039  \\
&(41 $\times$ 41 $\times$ 41)             &  4.9112E-5  &  4.019     \\
\hline\hline
 \end{tabular}
\end{table}
}
\subsection{Problem-2}
\begin{figure}[htbp]
 \centering
 \vspace*{5pt}%
 \hspace*{\fill}%
\begin{subfigure}{0.50\textwidth}     
    \centering
    \includegraphics[width=\textwidth]{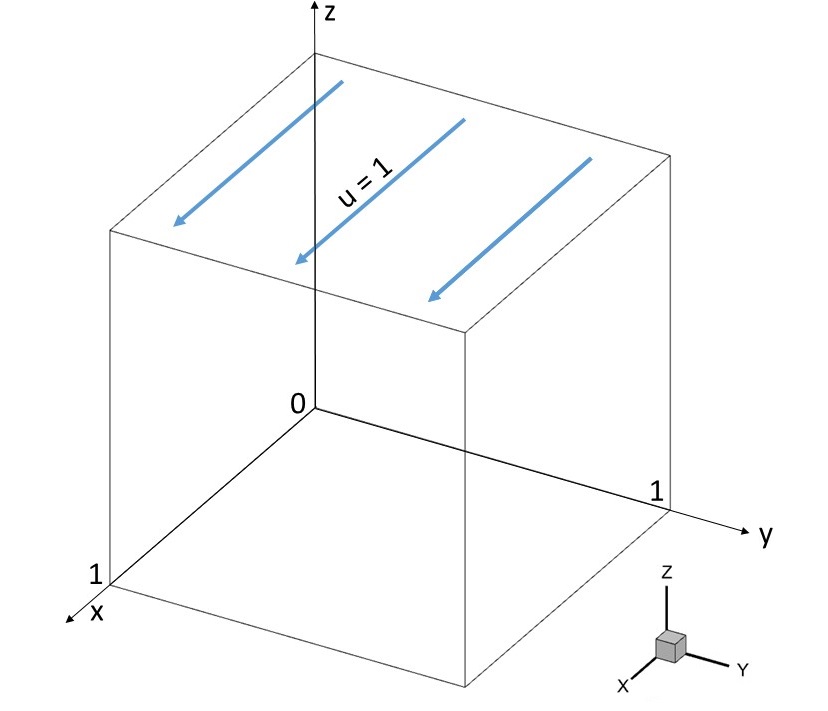}%
    \captionsetup{skip=5pt}%
    \caption{(a)}
    \label{fig:Cavity_3d}
  \end{subfigure}%
 \begin{subfigure}{0.50\textwidth}        
   \centering
    \includegraphics[width=\textwidth]{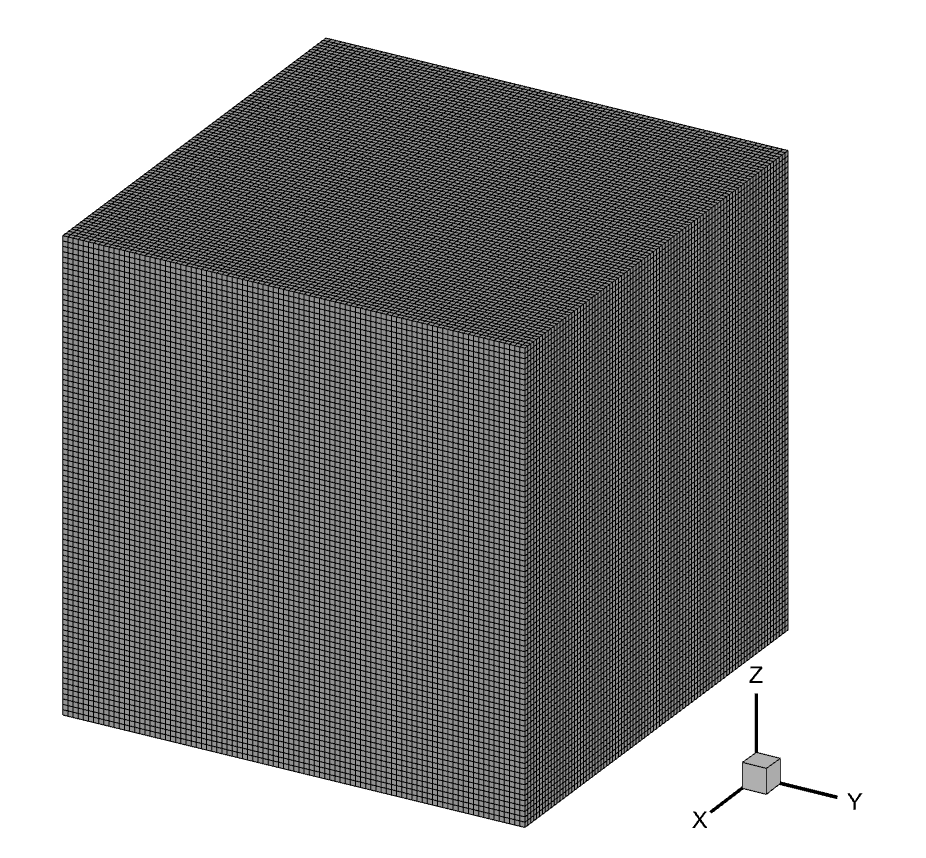}%
    \captionsetup{skip=5pt}%
    \caption{(b)}
    \label{fig:3D_grids_91_91_91}
  \end{subfigure}
  \hspace*{\fill}

  \vspace*{8pt}%
  \hspace*{\fill}%
  \caption{ (a) Illustration of the configuration in the 3D lid-driven cavity scenario and (b) View of the grids at a resolution of $91\times91\times91$.}
  \label{fig:Lid_cavity_grid}
\end{figure}

To check the ability of the new pressure correction technique, we have considered here the three-dimensional lid-driven cavity flow\cite{Ku_1987, Jiang_1994, Shu_2003, Zunic_2006, De_2009, Feldman_2010}. The scenario of a square cavity driven by the top lid has gained prominence as a widely accepted benchmark in the realm of fluid dynamics for evaluating the accuracy and performance of numerical methods used to solve the incompressible Navier-Stokes ($\mathrm{N}-\mathrm{S}$) equations. Its significance stems from the fact that this simple geometric configuration encapsulates various fluid mechanical phenomena typically observed in incompressible viscous flows.\\
In the 2D lid-driven square cavity problem, a square cavity is filled with an incompressible fluid, and the motion of the fluid is induced by imposing a constant velocity on the top boundary of the cavity (the top lid). The other three boundaries are stationary (no-slip condition), and the fluid is assumed to be viscous. Researchers often study various aspects of the lid-driven cavity problem, such as flow patterns, velocity profiles, pressure distributions, and circulation rates, to gain deeper insights into fluid dynamics and turbulence. In summary, the 2D lid-driven square cavity problem plays a crucial role in advancing the field of fluid mechanics and numerical methods. Its simplicity, universality, and ability to capture essential fluid phenomena make it an indispensable tool for researchers and engineers seeking to explore, develop, and validate numerical techniques for solving complex flow problems encountered in engineering, environmental sciences, and other relevant disciplines.\\
The 3D counterpart of the lid-driven square cavity problem presents a more complex and computationally demanding scenario compared to its 2D counterpart. The main challenge arises from the significant increase in the number of unknowns, leading to a substantial memory requirement and computational cost. Despite its complexity, the 3D lid-driven cubic cavity problem is of great interest in fluid dynamics and computational fluid dynamics (CFD) research due to its ability to capture a broader range of fluid phenomena, especially in three-dimensional flows. However, obtaining accurate numerical results for this problem is computationally intensive and often requires high-performance computing resources. In the 3D lid-driven cubic cavity problem, a cubic cavity is filled with an incompressible fluid. The motion of the fluid is induced by imposing a constant velocity to the top lid. The remaining five faces of the cube are stationary (no-slip condition). Figure \ref{fig:Lid_cavity_grid} illustrates a schematic diagram of the flow configuration within the cubic cavity. 
The initial velocity components in the entire domain are set to zero:
\[
u(x, y, z, t=0) = 0, \quad v(x, y, z, t=0) = 0, \quad w(x, y, z, t=0) = 0
\]
The boundary conditions for the 3D lid-driven cavity problem are as follows:\\
1. Lid-driven condition:
On the top face of the cubic cavity, a uniform velocity is prescribed along $x$ directions, creating a lid-driven flow. i.e:
\[
u(x, y, 1, t) = U_0, \quad v(x, y, 1, t) = 0, \quad w(x, y, 1, t) = 0
\]
Here, \( U_0 \) represents the velocity of the lid in the \( x \)-direction.\\
2. No-slip condition on five other faces:
\[
U(0, y, z, t) = 0, \quad U(1, y, z, t) = 0, \quad U(x, 0, z, t) = 0, \quad \] \\\[ U(x, 1, z, t) = 0, \quad U(x, y, 0, t) = 0
\]
Where, $U=(u,v,w)$. Neumann boundary conditions are applied for the pressure. 
\subsubsection{Grid Independence Study}
\begin{figure}[htbp]
 \centering
 \vspace*{5pt}%
 \hspace*{\fill}%
\begin{subfigure}{0.50\textwidth}     
    \centering
    \includegraphics[width=\textwidth]{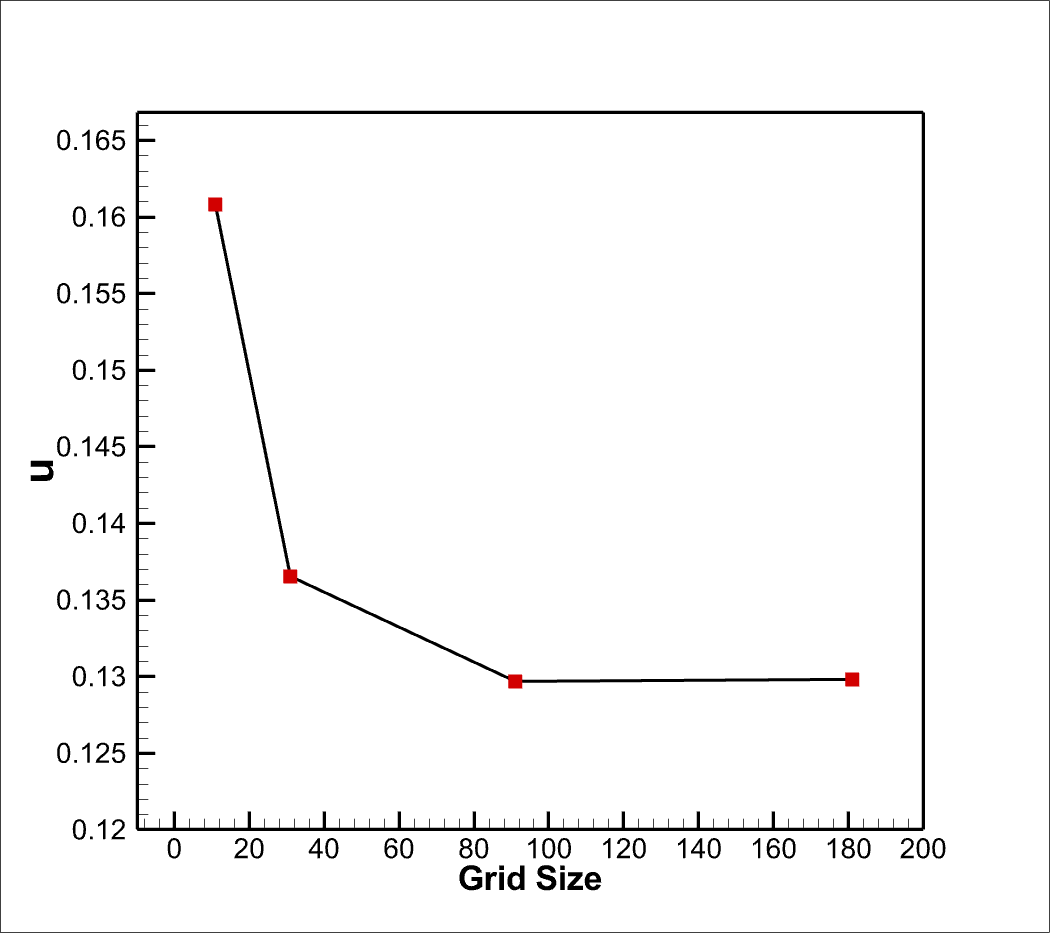}%
    \captionsetup{skip=5pt}%
    \caption{(a)}
    \label{fig:Grid_vs_Stream}
  \end{subfigure}%
 \begin{subfigure}{0.50\textwidth}        
   \centering
    \includegraphics[width=\textwidth]{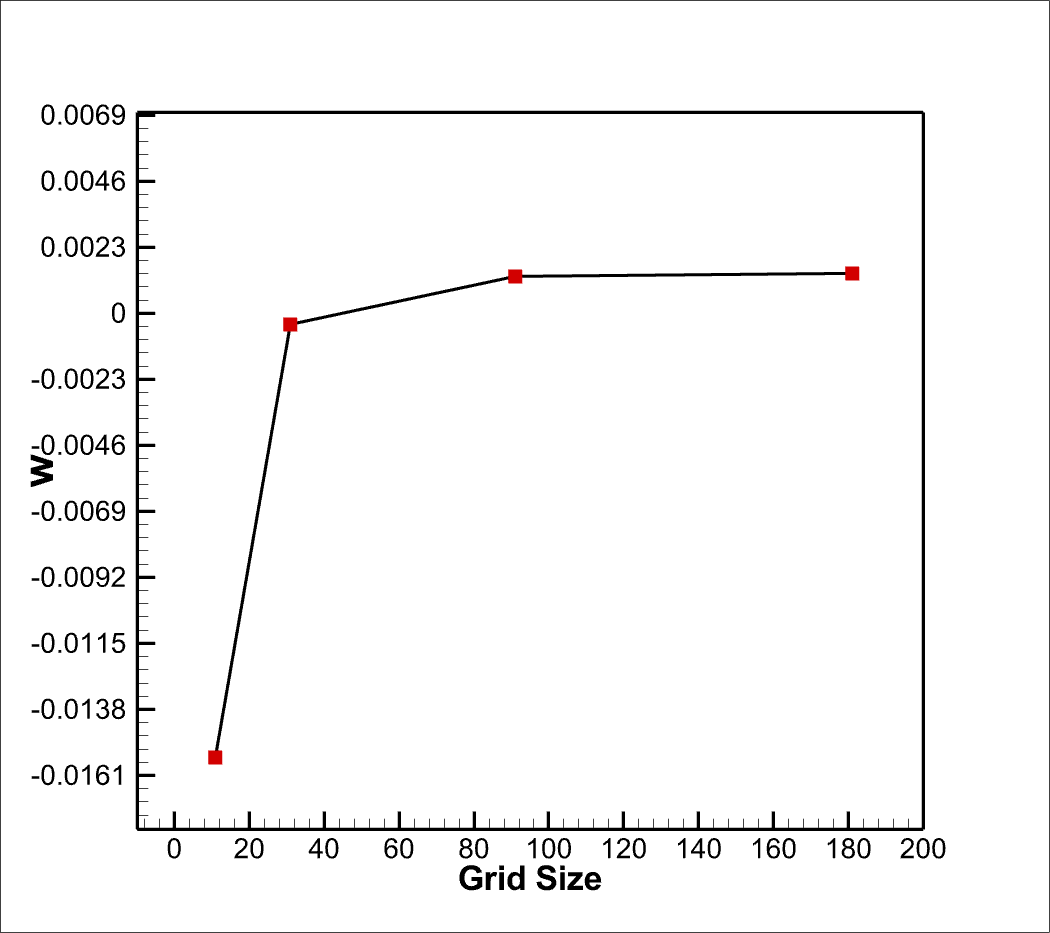}%
    \captionsetup{skip=5pt}%
    \caption{(b)}
    \label{fig:Pressure_vs_Grid}
  \end{subfigure}
  \hspace*{\fill}

  \vspace*{8pt}%
  \hspace*{\fill}%
  \caption{Velocity values along the different grid sizes at an observing point (0.75, 0.75, 0.75) (a) $u$ and (b) $w$}
  \label{fig:grid_ind}
\end{figure}

In order to guarantee the accuracy and reliability of our numerical simulations, we carry out an analysis of grid independence. The primary objective of this investigation is to assess how the outcomes respond to changes in grid resolutions, ultimately pinpointing the minimum grid size necessary for obtaining reliable solutions. The process of grid refinement analysis consists of progressively enhancing the grid resolution while maintaining consistent values for all other simulation parameters. Starting with an initial grid resolution, a series of grid refinements are performed by successively dividing the cell size in each dimension. The grid resolutions considered in this study are $11 \times 11 \times 11$, $31 \times 31 \times 31$, $91 \times 91 \times 91$, and $181 \times 181 \times 181$. 
For each grid resolution, simulations are conducted to compute the values of velocities($u,v,w$) at an observation point (0.75, 0.75, 0.75). The calculated values are presented in Table-\ref{grid_independent_test} and visualized in Figures \ref{fig:Grid_vs_Stream} and \ref{fig:Pressure_vs_Grid}. To ensure consistency, the same physical problem setup, initial conditions, and boundary conditions are maintained throughout the grid refinement process.  
{\small\begin{table}[htbp]
\caption{\small Velocities values at a point of observation $(0.75, 0.75, 0.75)$ at time $=15$ for Reynold number $Re = 100$, and  $\Delta$ t = 0.02 by employing four distinct sizes of grid}\label{grid_independent_test}
\centering
 \begin{tabular}{ccccc}  \hline \hline
M$\times$N$\times$O    &    $u$   &  $v$  &  $w$ & Max. Relative Error ($\%$)    \\ \hline
(11 $\times$ 11 $\times$ 11)           &  0.16080 &   -0.00002    &   -0.01547  & 15.07$\%$ \\
(31 $\times$ 31 $\times$ 31)         &  0.13656  &  -0.00203  &   -0.00038  & 5.02$\%$\\
(91 $\times$ 91 $\times$ 91)              &  0.12970 &  -0.00298  &   0.00128  & 0.29$\%$  \\
(181 $\times$ 181 $\times$ 181)              &  0.12981 &  -0.00305 &   0.00137 & --   \\
\hline
 \end{tabular}
\end{table}
}

Figures \ref{fig:Grid_vs_Stream} and \ref{fig:Pressure_vs_Grid} depict the changes in $u$- and $w$-velocities concerning the grid resolution. The outcomes exhibit a clear trend of convergence as the grid becomes more refined, signifying a decrease in solution errors as the grid resolution increases. Table \ref{grid_independent_test} provides a quantitative assessment of grid convergence, presenting computed velocity values at an observation point for each grid resolution. The values demonstrate the reduction in differences among successive grid resolutions and provide an indication of the point of grid independence. Table \ref{grid_independent_test} reveals that the maximum relative error between the last two grid sizes is just 0.29$\%$. It's evident that a grid size of (91 $\times$ 91 $\times$ 91) 
 is sufficient to accurately capture the flow phenomena. As a result, we adopt a grid size of (91 $\times$ 91 $\times$ 91) for all subsequent calculations.
\begin{figure}[htbp]
 \centering
 \vspace*{5pt}%
 \hspace*{\fill}%
\begin{subfigure}{0.50\textwidth}     
    \centering
    \includegraphics[width=\textwidth]{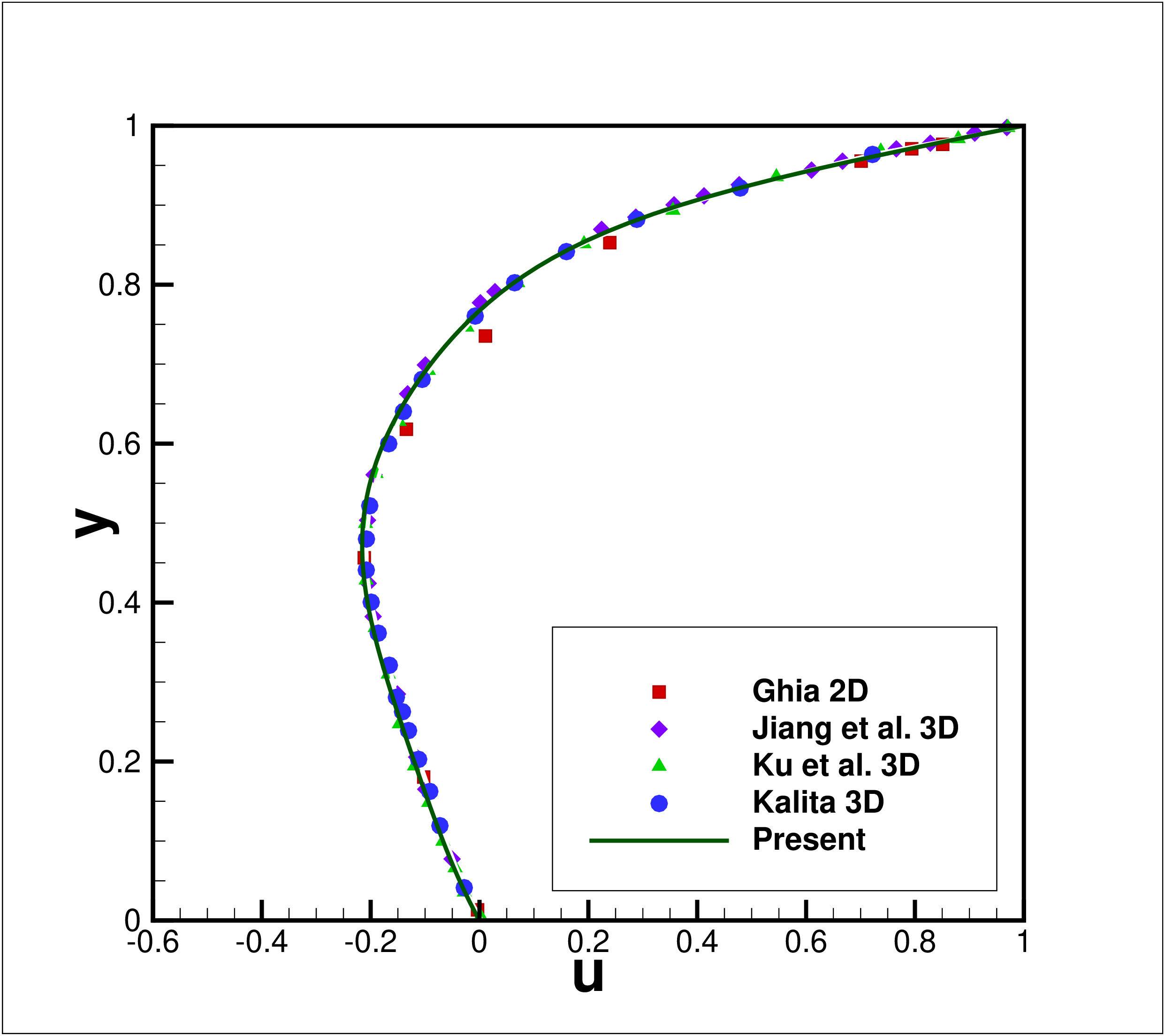}%
    \captionsetup{skip=5pt}%
    \caption{(a)}
    \label{fig:Final_comparison_u_y}
  \end{subfigure}%
 \begin{subfigure}{0.50\textwidth}        
   \centering
    \includegraphics[width=\textwidth]{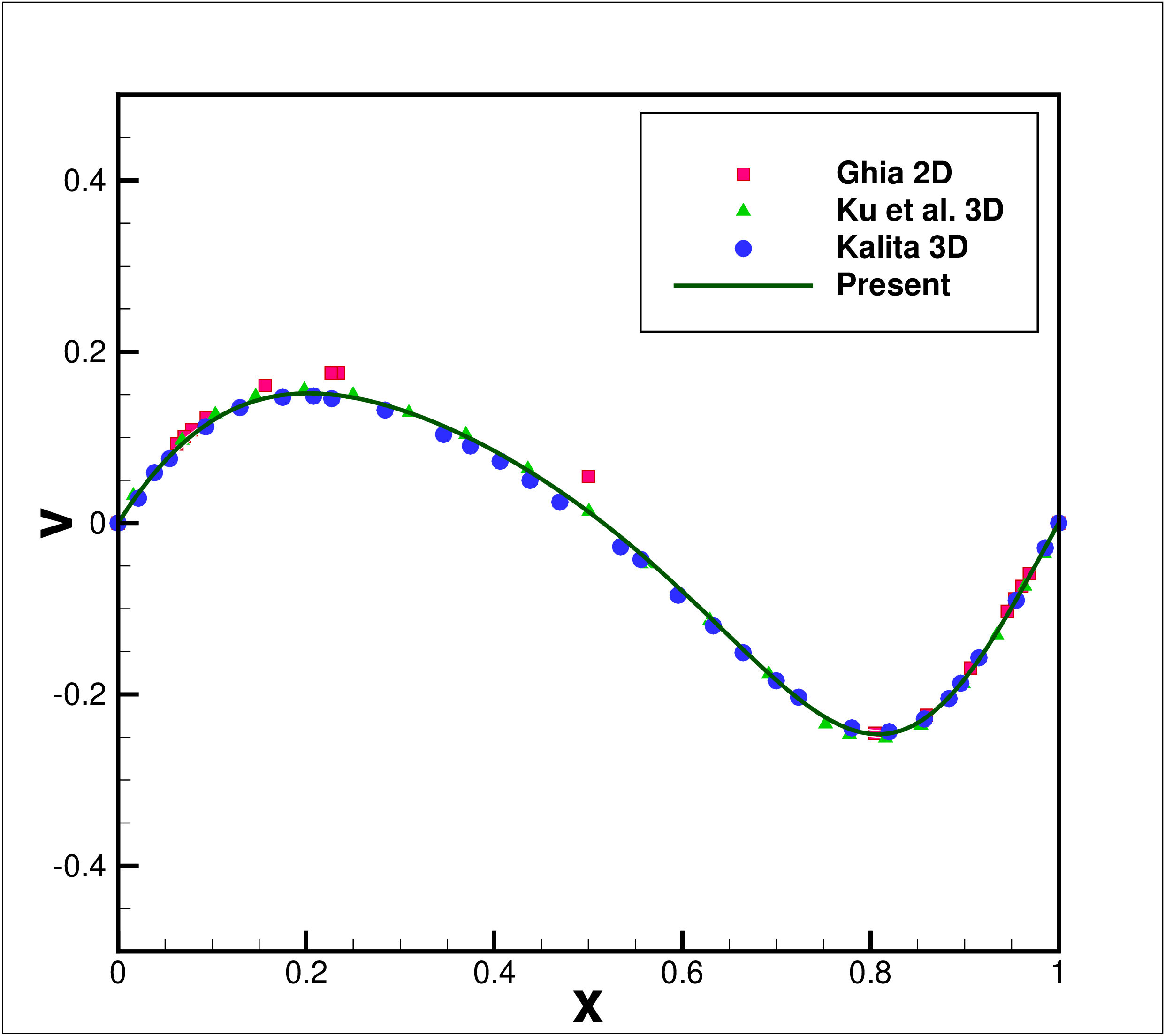}%
    \captionsetup{skip=5pt}%
    \caption{(b)}
    \label{fig:Final_comparison_x_v}
  \end{subfigure}
  \hspace*{\fill}

  \vspace*{8pt}%
  \hspace*{\fill}%
  \caption{ Comparison between the 2D and 3D simulations\cite{Ghia_1982,Ku_1987, Jiang_1994,Kalita_2014} at $Re=100$ (a) Horizontal Velocity Profile along the Vertical Centerline (b)  Vertical Velocity Profile along the Horizontal Centerline}
  \label{fig:velocity_comparison}
\end{figure}

\subsubsection{Results and Analysis}
After successfully completing the grid independence test and determining the optimal grid size, we proceeded to compute the steady-state solutions for the three-dimensional lid-driven cavity flow problem. The steady state is considered to be reached when the following inequality is satisfied:\\
$\max \sqrt{\left(\Delta u_{i j k}\right)^2+\left(\Delta v_{i j k}\right)^2+\left(\Delta w_{i j k}\right)^2}< 10^{-6}$\\
Here, $\Delta u_{i j k} = u_{i j k}^{(n+1)}-u_{i j k}^{(n)}$, etc. with $(n)$ denotes the previous time level and $(n+1)$ represents the current time level. Our simulations encompassed Reynolds numbers of $100, 400,$ and $1000$. The numerical simulations are carried out using the higher order super compact method, incorporating the novel effective pressure correction technique. This new pressure correction technique not only reduces computational costs during each time increment but also provides accurate and robust results. 
\begin{figure}[htbp]
 \centering
    \includegraphics[width=\textwidth]{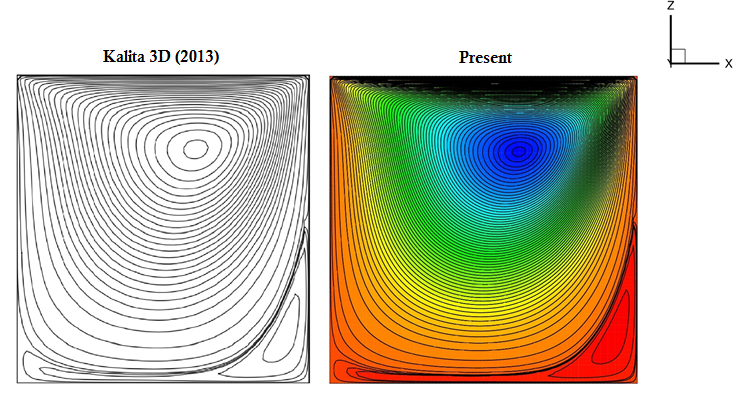}%
  \hspace*{\fill}%
  \caption{Comparison of Streamlines plot at $y=0.5$ for $Re=100$}
  \label{fig:Re_100_Streamlines}
\end{figure}
\begin{figure}[htbp]
 \centering
    \includegraphics[width=\textwidth]{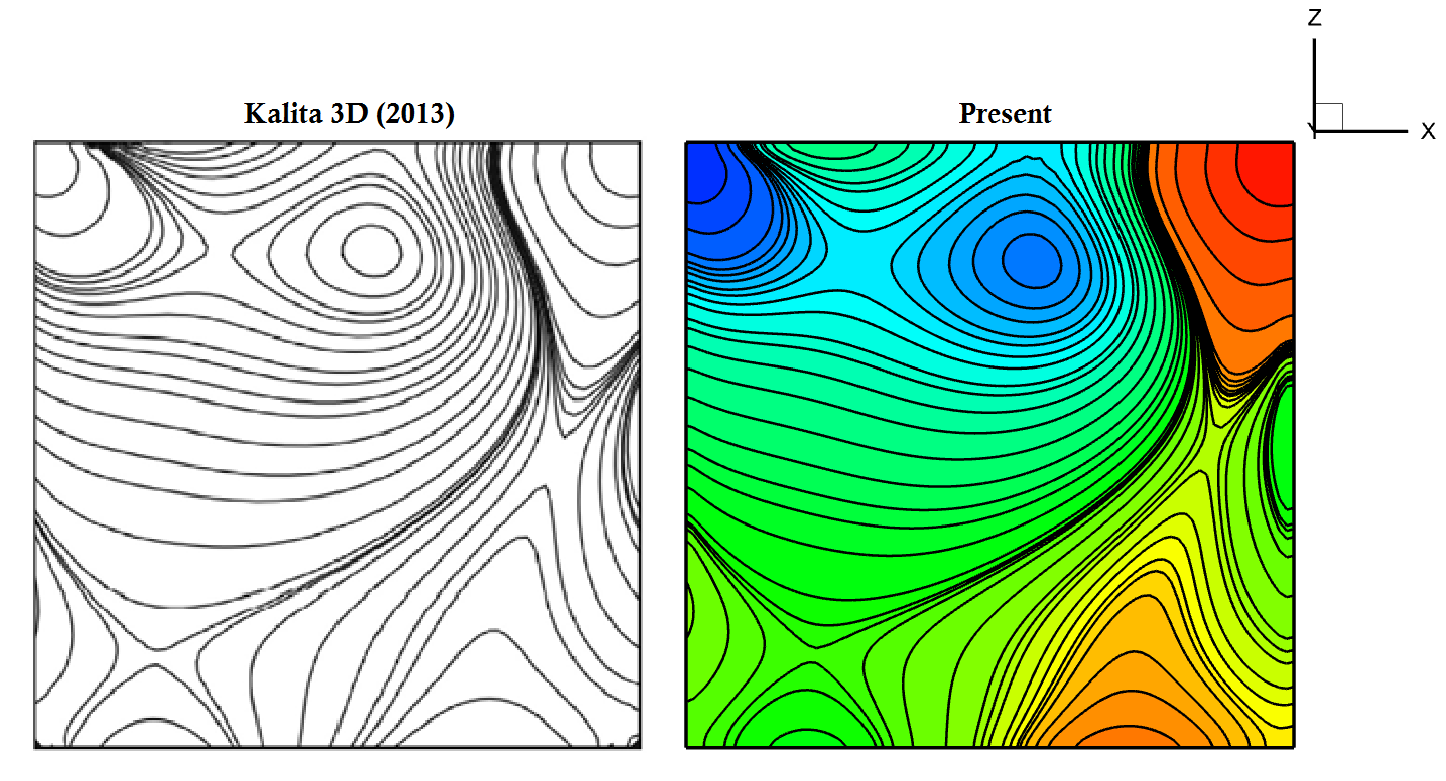}%
  \hspace*{\fill}%
  \caption{ Comparison of Pressure contours at $y=0.5$ for $Re=100$}
  \label{fig:Re_100_Pressure1}
\end{figure}
\begin{figure}[htbp]
 \centering
 \vspace*{5pt}%
 \hspace*{\fill}%
\begin{subfigure}{0.33\textwidth}     
    \centering
    \includegraphics[width=\textwidth]{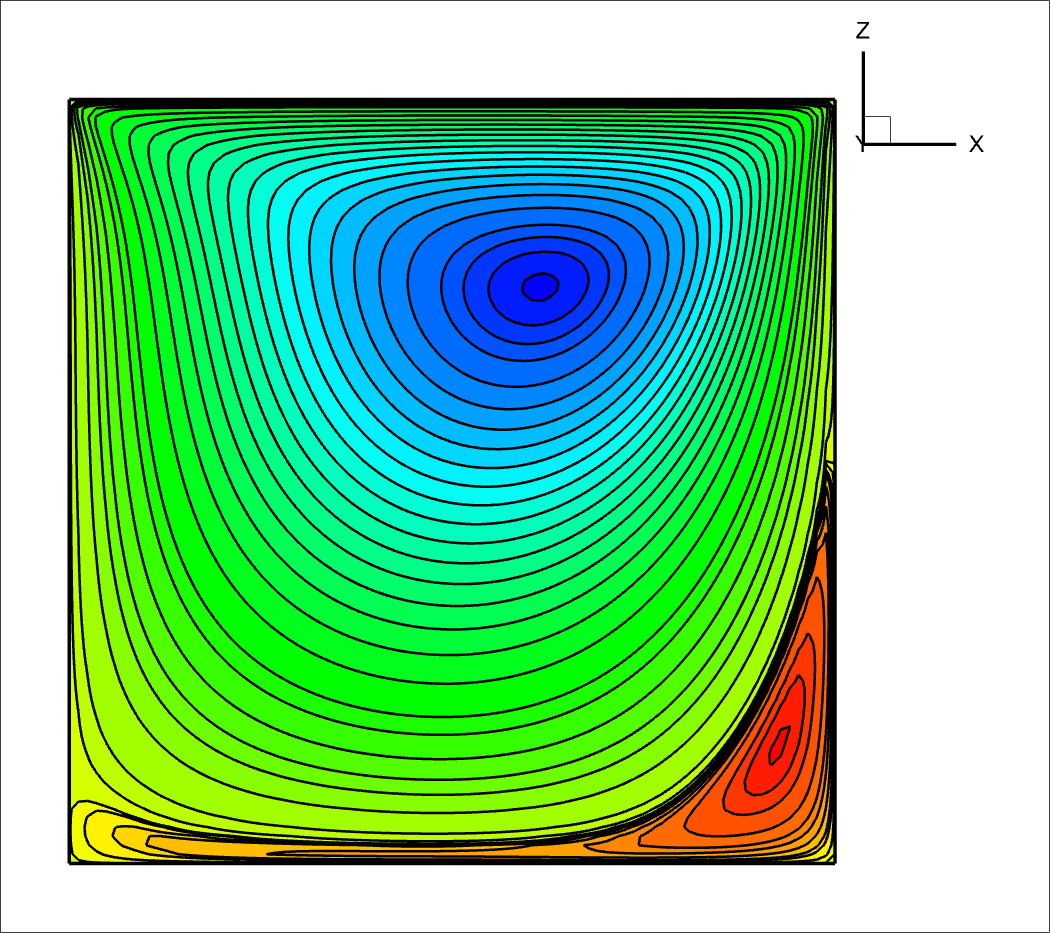}%
    \captionsetup{skip=5pt}%
    \caption{(a)}
    \label{fig:Streamlines_Re_100}
  \end{subfigure}%
 \begin{subfigure}{0.33\textwidth}        
   \centering
    \includegraphics[width=\textwidth]{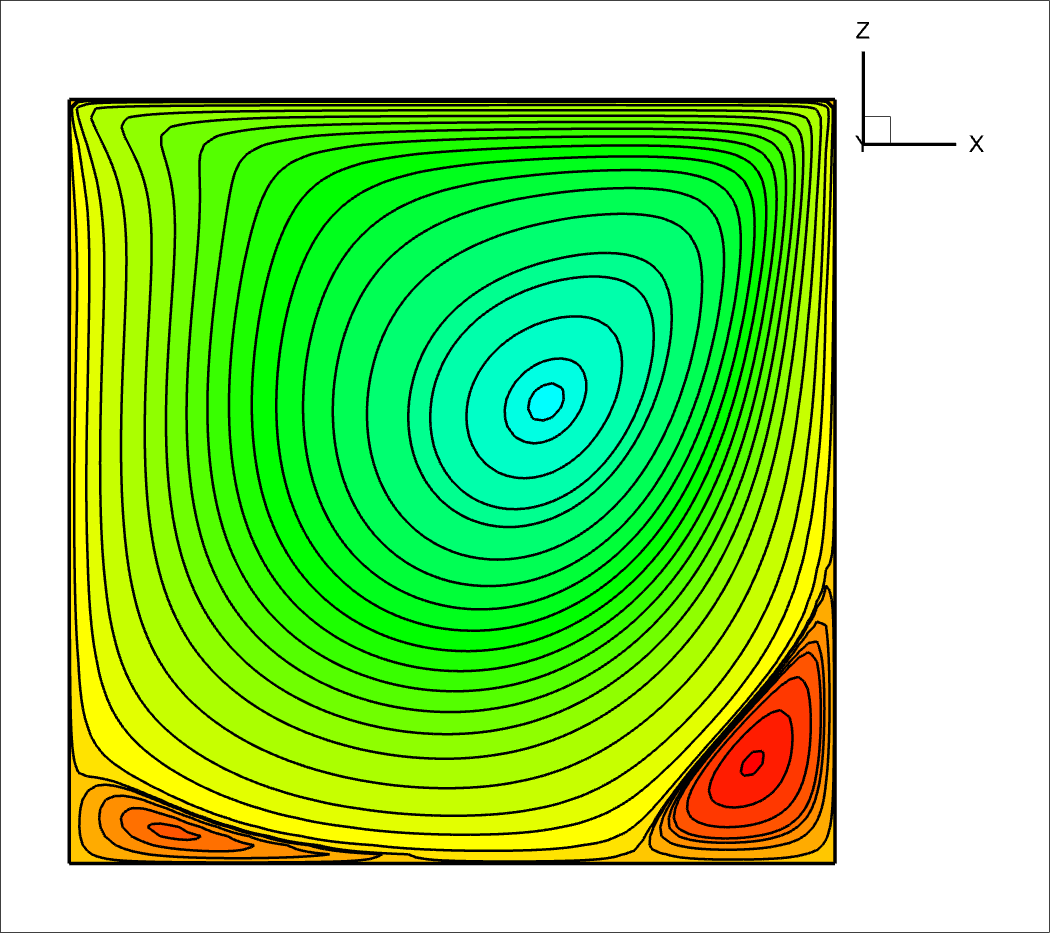}%
    \captionsetup{skip=5pt}%
    \caption{(b)}
    \label{fig:Streamlines_Re_400}
  \end{subfigure}
   \begin{subfigure}{0.33\textwidth}        
   \centering
    \includegraphics[width=\textwidth]{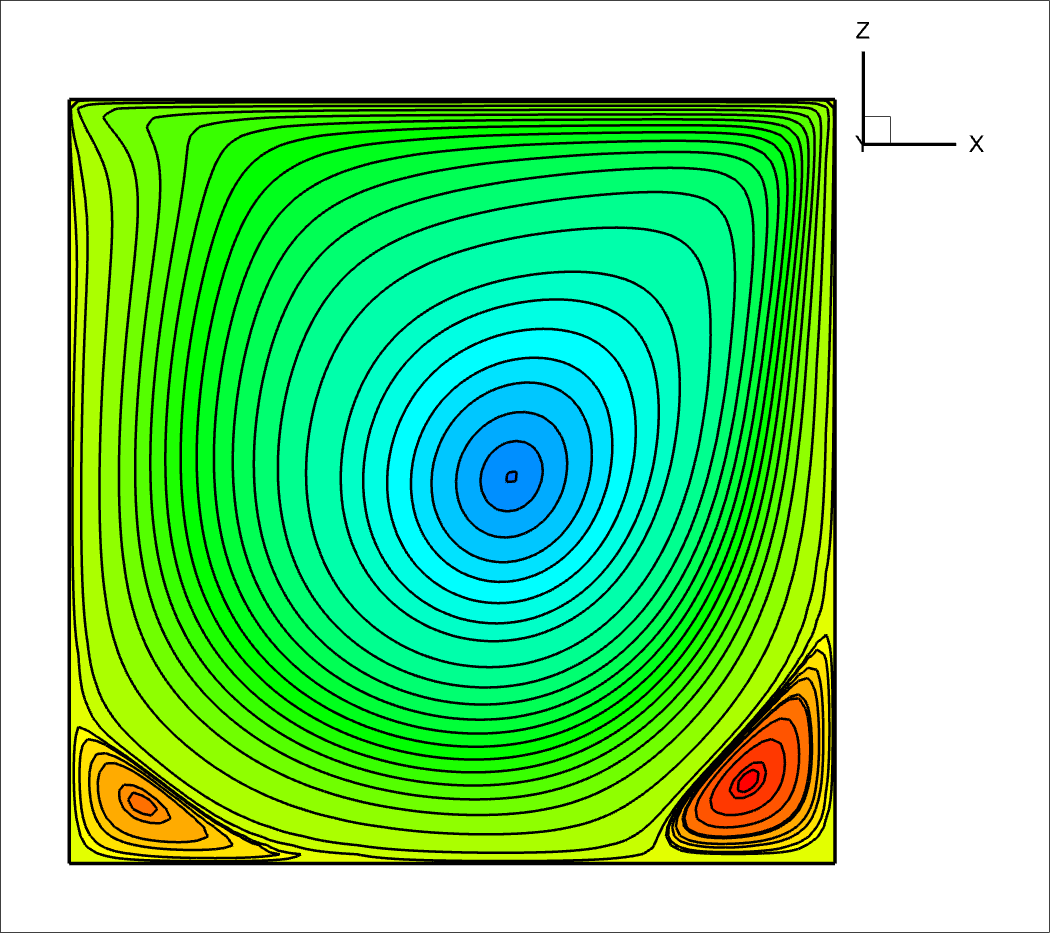}%
    \captionsetup{skip=5pt}%
    \caption{(c)}
    \label{fig:Streamlines_Re_1000}
  \end{subfigure}
  \hspace*{\fill}

  \vspace*{8pt}%
  \hspace*{\fill}%
  \caption{Streamline visualization on the $y = 0.5$ plane for the lid-driven cubic cavity problem: (a) $Re=100$, (b) $Re=400$, (c) $Re=1000$}
  \label{fig:Streamlines_lid_driven}
\end{figure}

{\small\begin{table}[htbp]
\caption{\small Comparison of the center of the primary vortex for different Reynold number}\label{Primary_Vortex}
\centering
 \begin{tabular}{ccc}  \hline  \hline
& $Re$    &   Primary Vortex       \\ \hline 
&100          &  (0.616, 0.5, 0.755) \\
&400         &  (0.622, 0.5, 0.606)  \\
&1000             &  (0.579, 0.5, 0.508)      \\
\hline\hline
 \end{tabular}
\end{table}
}
\begin{figure}[htbp]
 \centering
 \vspace*{5pt}%
 \hspace*{\fill}%
\begin{subfigure}{0.33\textwidth}     
    \centering
    \includegraphics[width=\textwidth]{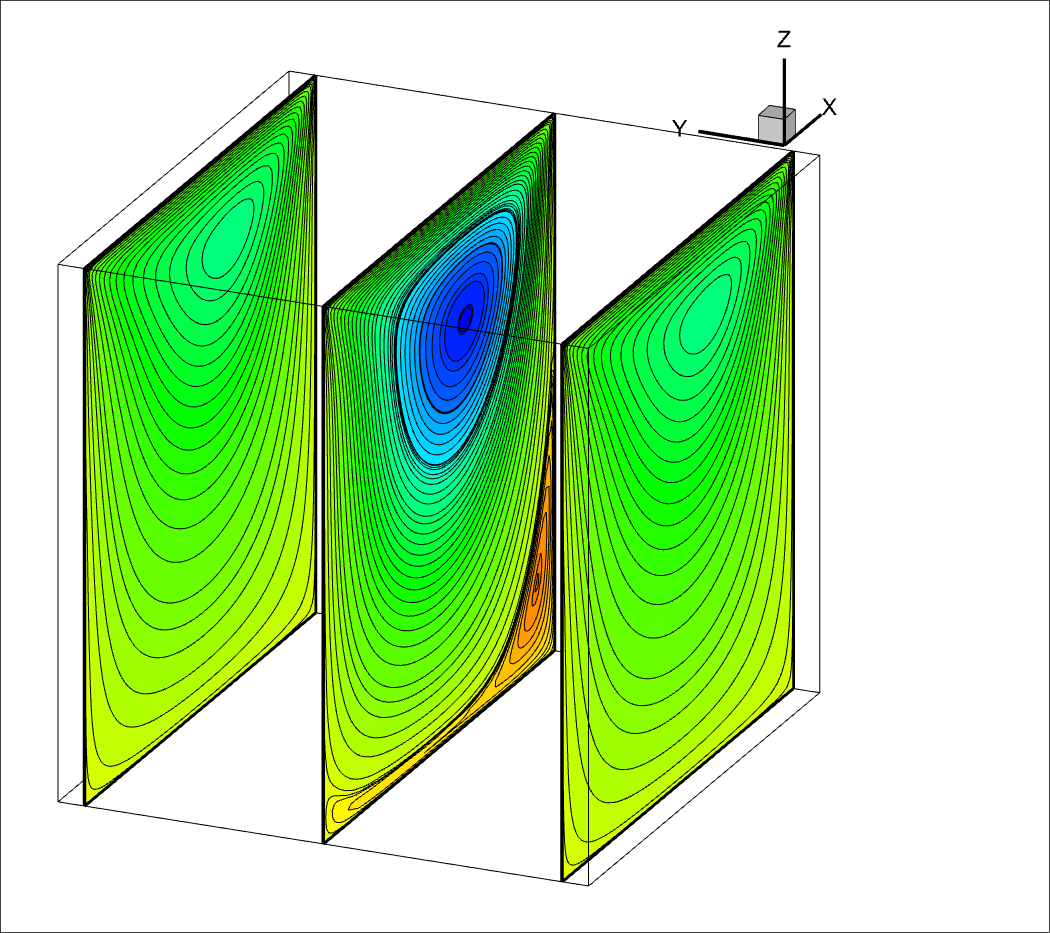}%
    \captionsetup{skip=5pt}%
    \caption{(a)}
    \label{fig:Streamlines_Re_100_3_Plane}
  \end{subfigure}%
 \begin{subfigure}{0.33\textwidth}        
   \centering
    \includegraphics[width=\textwidth]{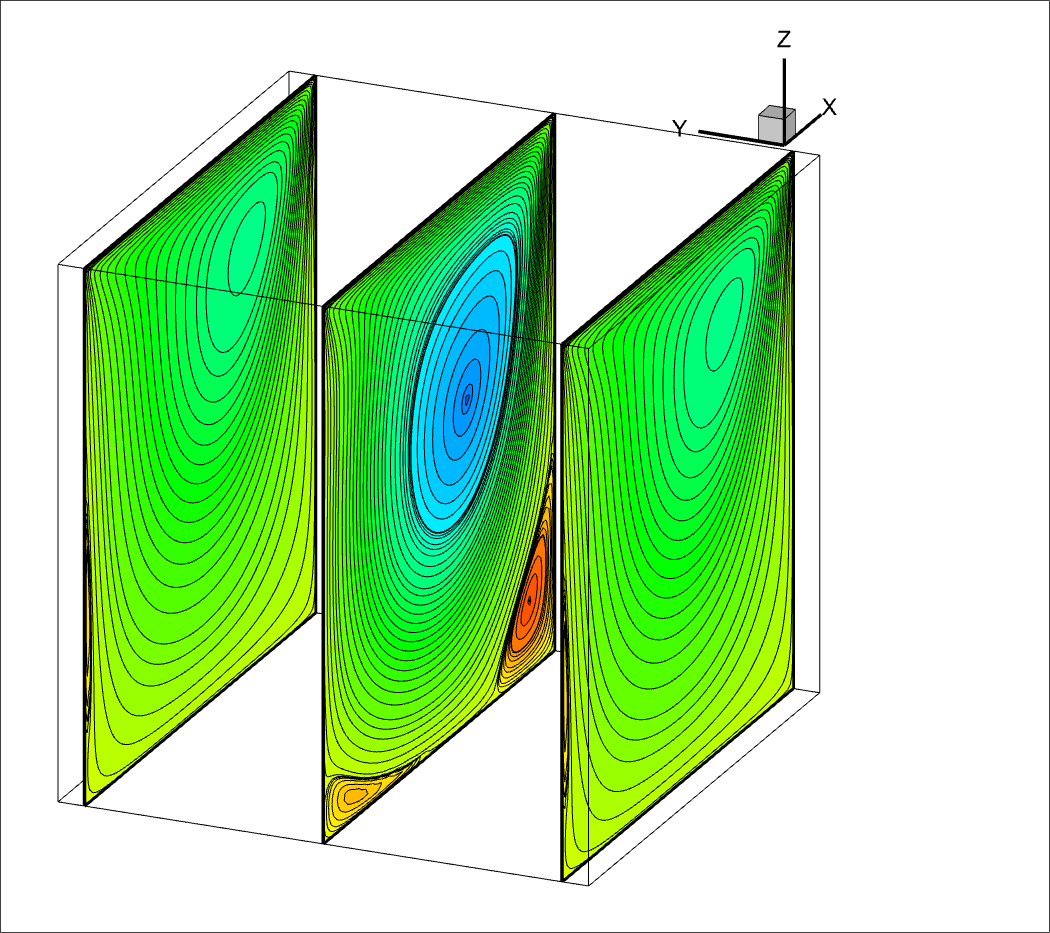}%
    \captionsetup{skip=5pt}%
    \caption{(b)}
    \label{fig:Streamlines_Re_400_3_plane}
  \end{subfigure}
   \begin{subfigure}{0.33\textwidth}        
   \centering
    \includegraphics[width=\textwidth]{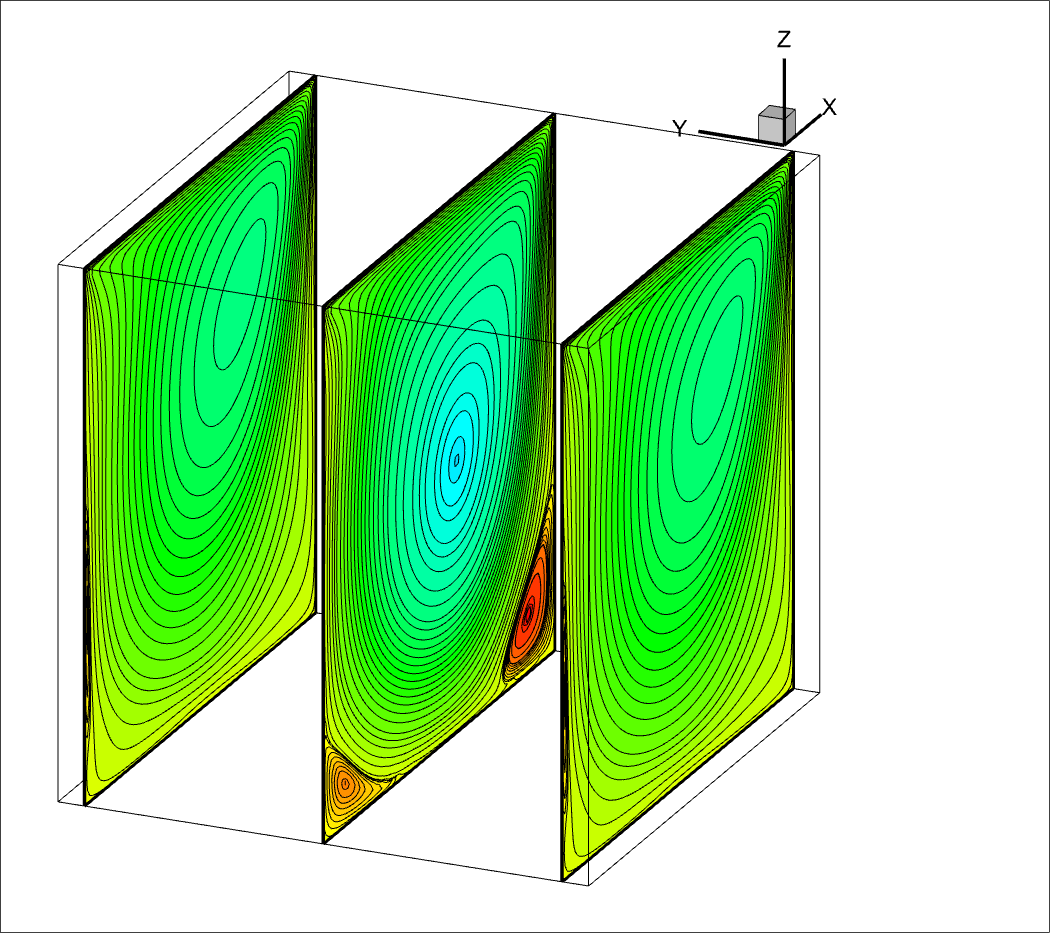}%
    \captionsetup{skip=5pt}%
    \caption{(c)}
    \label{fig:Streamlines_Re_1000_3_plane}
  \end{subfigure}
  \hspace*{\fill}

  \vspace*{8pt}%
  \hspace*{\fill}%
  \caption{Streamline visualization on the three slices ($y = 0.05$, $y = 0.5$ and $y = 0.95$ planes) for the lid-driven cubic cavity problem: (a) $Re=100$, (b) $Re=400$, (c) $Re=1000$}
  \label{fig:Streamlines_lid_driven_3_plane}
\end{figure}

Figure \ref{fig:velocity_comparison} presents a detailed comparison between our computational 3D results and existing data (2D and 3D) \cite{Ghia_1982,Ku_1987, Jiang_1994,Kalita_2014} for the variation of velocity components, specifically $u$ vs. $y$ and $x$ vs. $v$, in the lid-driven cavity problem. The plot illustrates the excellent agreement between our numerical simulations and the pre-existing data. The close match between the sets of results indicates the accuracy and reliability of our computational approach in capturing the intricate velocity variations within the cavity. In Figures \ref{fig:Re_100_Streamlines} and \ref{fig:Re_100_Pressure1}, we present the comparison of the streamlines and pressure plots obtained from our simulations for $y=0.5$ and $Re=100$ with the results reported by Kalita (at 65$\times$ 65 $\times$ 65 grid size)\cite{Kalita_2014}. The plots clearly demonstrate excellent agreement between our computed results and those obtained by Kalita \cite{Kalita_2014}. The streamlines and pressure distributions show very similar results, which indicates that our numerical approach accurately captured the flow behavior in the lid-driven cavity.\\

\begin{figure}[htbp]
 \centering
 \vspace*{5pt}%
 \hspace*{\fill}%
\begin{subfigure}{0.33\textwidth}     
    \centering
    \includegraphics[width=\textwidth]{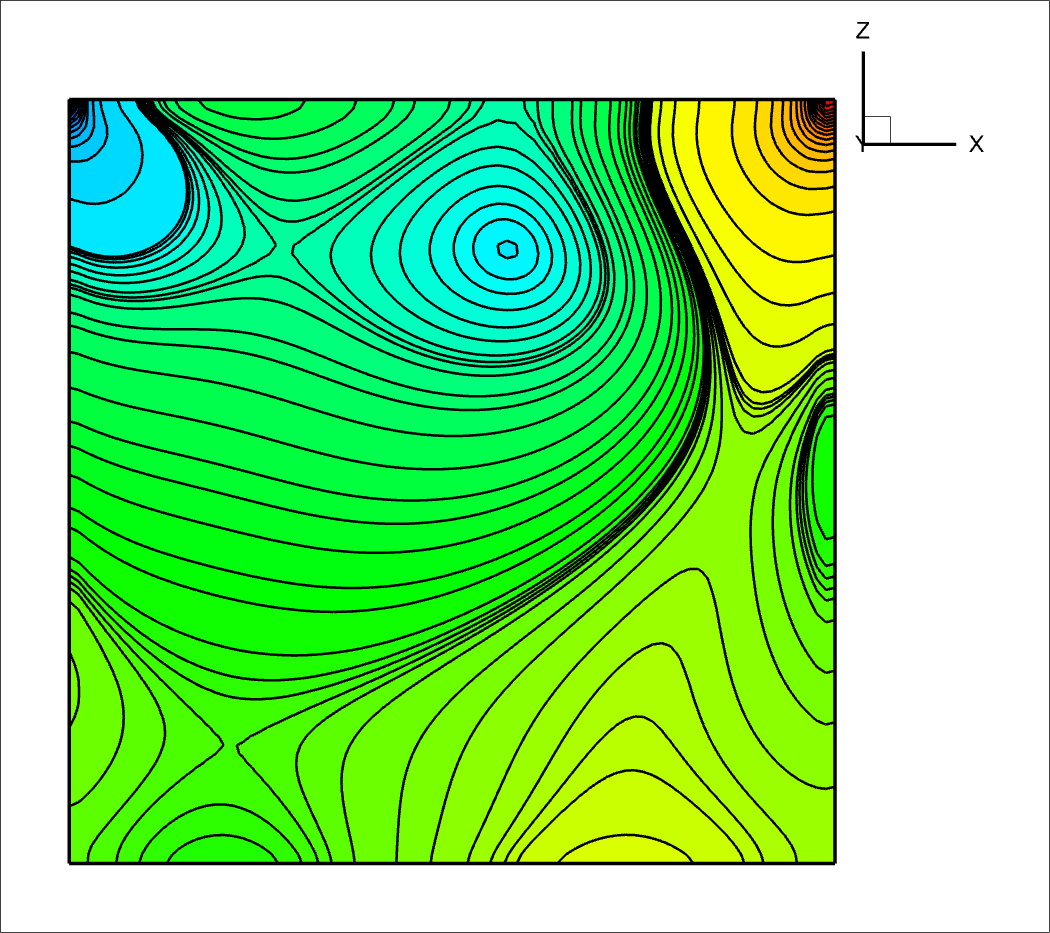}%
    \captionsetup{skip=5pt}%
    \caption{(a)}
    \label{fig:Pressure_Re_100}
  \end{subfigure}%
 \begin{subfigure}{0.33\textwidth}        
   \centering
    \includegraphics[width=\textwidth]{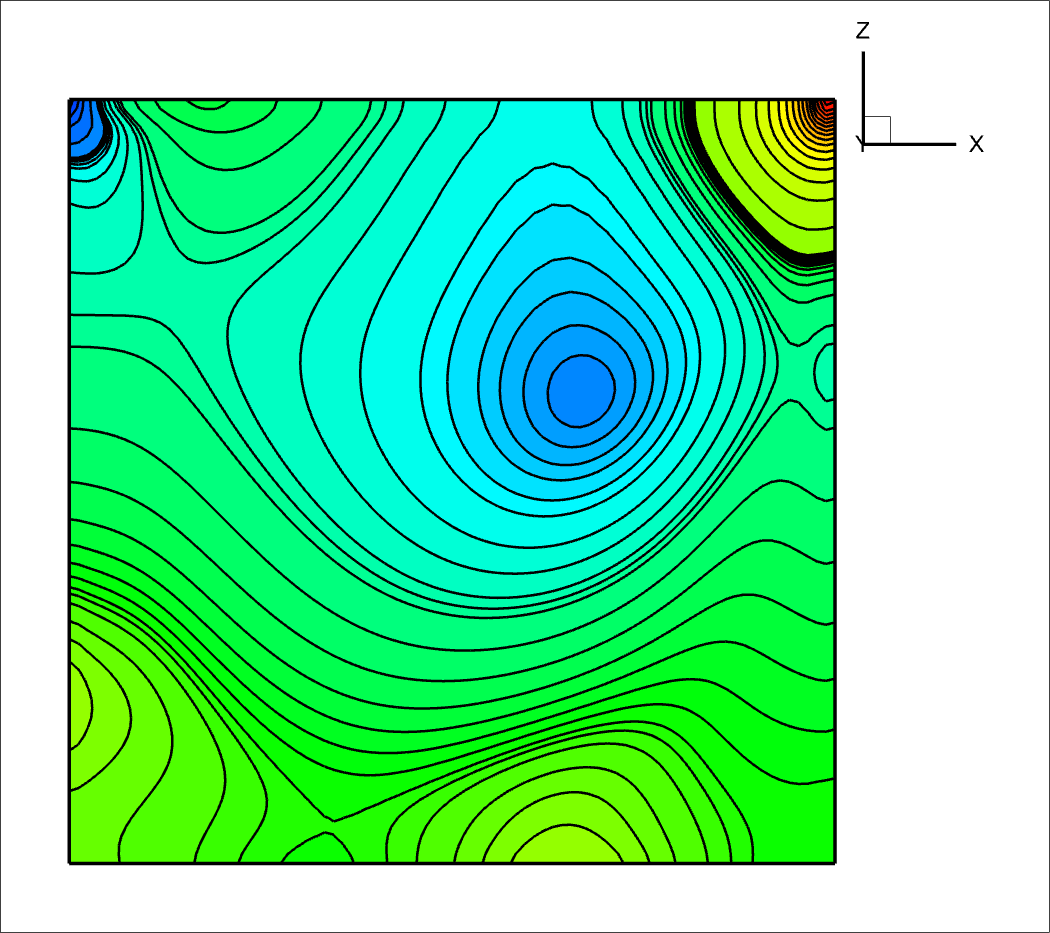}%
    \captionsetup{skip=5pt}%
    \caption{(b)}
    \label{fig:Pressure_Re_400}
  \end{subfigure}
   \begin{subfigure}{0.33\textwidth}        
   \centering
    \includegraphics[width=\textwidth]{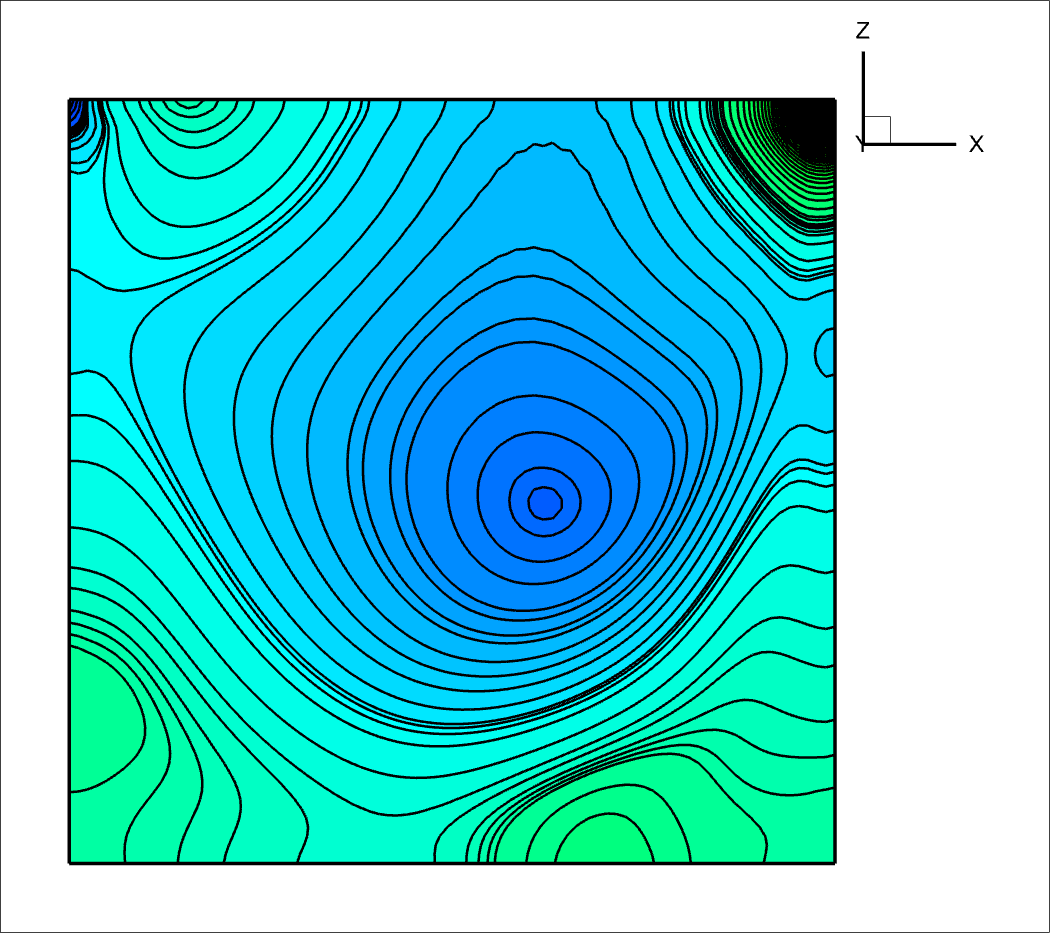}%
    \captionsetup{skip=5pt}%
    \caption{(c)}
    \label{fig:Pressure_Re_1000}
  \end{subfigure}
  \hspace*{\fill}

  \vspace*{8pt}%
  \hspace*{\fill}%
  \caption{Pressure distribution on the $y = 0.5$ plane for the lid-driven cubic cavity problem: (a) $Re=100$, (b) $Re=400$, (c) $Re=1000$ }
  \label{fig:Pressure_lid_driven}
\end{figure}


{\small\begin{table}[htbp]
\caption{\small Comparative Analysis of Various Parameters in the 3D Lid-Driven Cavity Problem-2}\label{Pressure_Iterations}
\centering
 \begin{tabular}{cccccc}  \hline \hline
 &    &    Without Using Proposed Pressure correction Technique &    &    & \\ \hline
& $Re$    &   $\lambda$    &  $\Delta$t  &  Iterations  &    \\ \hline 
&10          &  0.1  &    0.01  &    12\\
&100         &  0.06   & 0.005  &    20\\
&400             &  0.02  &  0.001 &    30    \\
&1000             &  0.02  &  0.000125 &    50    \\\hline 
&    &     Using Proposed Pressure correction Technique &    &    &    \\ \hline 
&10          &  0.1  &    0.02  &    2\\
&100         &  0.06   & 0.005  &    3\\
&400             &  0.02  &  0.002 &    4    \\
&1000             &  0.02  &  0.001 &    4   \\
\hline\hline
 \end{tabular}
\end{table}
}

Figure \ref{fig:Streamlines_lid_driven} showcases streamline contours along the $y = 0.5$ plane, while Figure \ref{fig:Streamlines_lid_driven_3_plane} exhibits streamline contours at three distinct plane slices. This visualization offers insight across an extensive spectrum of Reynolds numbers ($Re = 100$, $400$ and $1,000$). For all the considered cases, a prominent secondary vortex emerges at the bottom right corner of the cavity.  Notably, the emergence of the secondary vortex in the lower left corner becomes evident only when $Re$ reaches or surpasses the value of $400$. Comparing the size and shape of these secondary vortices with those observed in the two-dimensional cavity, we note the differences in the sizes and shapes of the vortices. These illustrations also serve to highlight that as the Reynolds number (Re) escalates, there is a discernible tendency for the center of the primary vortex to migrate closer to the geometric center. This trend is readily discernible through the data presented in Table \ref{Primary_Vortex}. This behavior is consistent with the expected influence of higher $Re$ on the flow dynamics, leading to changes in the flow patterns and vortex behavior. 

In Figure \ref{fig:Pressure_lid_driven}, we show the pressure contours in the middle of the $y$-plane ($y=0.5$) for different Reynolds numbers ($Re$ = 100, 400, and 1,000). The pressure contours offer valuable insights into the pressure distribution within the three-dimensional lid-driven cavity flow for varying Reynolds numbers. As evident from the plots, the pressure distribution exhibits intriguing variations, and the differences become more pronounced as the Reynolds number increases. The pressure contours provide crucial information about the flow behavior and allow us to discern the influence of the Reynolds number on the pressure distribution within the cavity. Overall, these findings shed light on the intricate flow phenomena occurring in the three-dimensional lid-driven cavity, providing valuable insights into how the flow structures evolve with varying Reynolds numbers. These observations contribute to a deeper understanding of the complex fluid dynamics associated with this problem and can serve as a basis for further investigations into more advanced turbulent flow scenarios.\\
We reproduced the results and noted a comparable number of pressure iterations to those documented in reference \cite{Kalita_2014} when we abstained from applying our proposed technique. Table \ref{Pressure_Iterations} presents a comprehensive comparison of the computational parameters used for solving the lid-driven cavity problem, both with and without the application of our proposed technique at $91 \times 91 \times 91$ grid size. The table shows the values of $\Delta t$ and the number of iterations required for pressure calculations for different Reynolds numbers. Pressure iterations are terminated when $|\nabla \cdot \mathbf{v}|_{\max}$  reaches a tolerance limit of $10^{-3}$. Remarkably, the results indicate a significant reduction in the number of iterations when our proposed technique is employed. This reduction in iterations directly translates to a substantial reduction in computational cost, making the proposed technique highly efficient in solving the problem. It's worth noting that for a Reynolds number of 1000, the conventional approach necessitates a time step size of $\Delta t = 0.000125$ for the convergence of the iterative method, whereas with the utilization of our pressure correction technique, the required time step is only $\Delta$ t = 0.001. At the same time, our technique takes only 4 inner iterations to converge to the desired accuracy in each time step, whereas the conventional approach needs 50 inner iterations for the same task. This also demonstrates the efficiency and computational cost-effectiveness of the proposed pressure correction technique. The improved computational efficiency achieved with our approach holds great promise for practical engineering and scientific simulations, as it allows for faster and more economical simulations while maintaining high accuracy and reliability.

\begin{figure}[htbp]
 \centering
 \vspace*{5pt}%
 \hspace*{\fill}%
\begin{subfigure}{0.50\textwidth}     
    \centering
    \includegraphics[width=\textwidth]{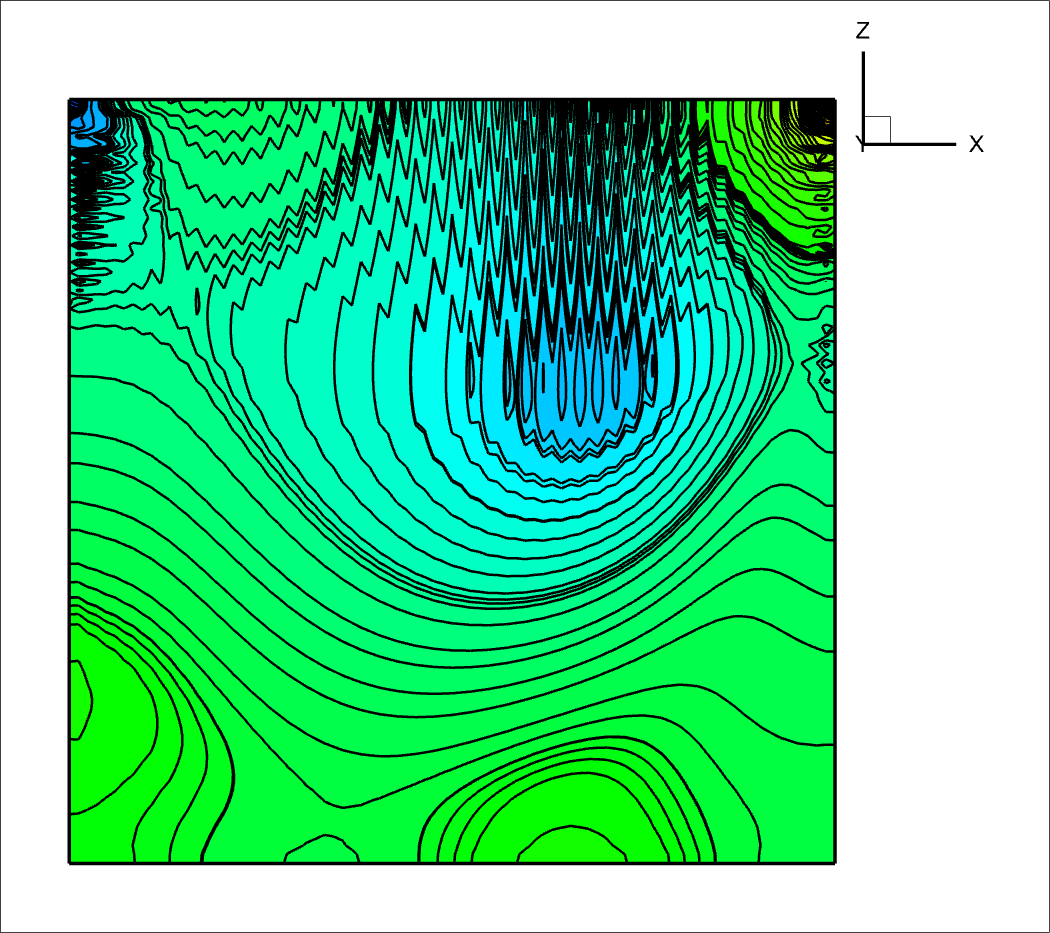}%
    \captionsetup{skip=5pt}%
    \caption{(a)}
    \label{fig:Pressure_Re_400_zig_zag}
  \end{subfigure}%
 \begin{subfigure}{0.50\textwidth}        
   \centering
    \includegraphics[width=\textwidth]{Pressure_Re_400.png}%
    \captionsetup{skip=5pt}%
    \caption{(b)}
    \label{fig:Pressure_Re_400c}
  \end{subfigure}
  \hspace*{\fill}

  \vspace*{8pt}%
  \hspace*{\fill}%
  \begin{subfigure}{0.50\textwidth}     
    \centering
    \includegraphics[width=\textwidth]{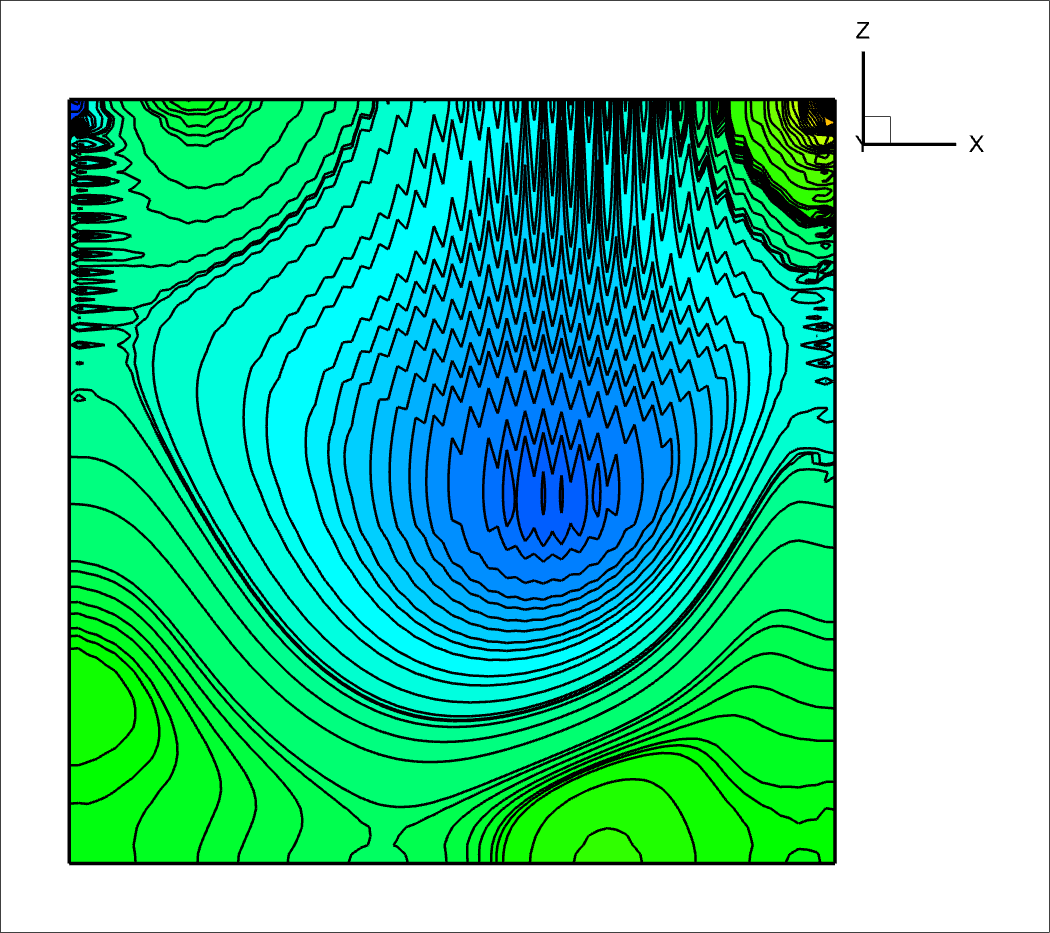}%
    \captionsetup{skip=5pt}%
    \caption{(c)}
    \label{fig:Pressure_Re_1000_zig_zag}
  \end{subfigure}%
  \begin{subfigure}{0.50\textwidth}     
    \centering
    \includegraphics[width=\textwidth]{Pressure_Re_1000.png}%
    \captionsetup{skip=5pt}%
    \caption{(d)}
    \label{fig:Pressure_Re_1000c}
  \end{subfigure}%
  \caption{Comparison of the pressure contours (a,c) without using (left) and (b,d) with using (right) our proposed pressure correction technique at the same number of pressure iterations for (a,b) $Re=400$, (c,d) $Re=1000$.}
  \label{fig:Pressure_zig_zag}
\end{figure}
Figure \ref{fig:Pressure_zig_zag} displays the pressure contours obtained for Reynolds numbers 400 and 1000 using the proposed Gaussian-based technique and the conventional method, both employing the same number of iterations. A keen observation reveals striking differences in the results between the two approaches. When not using our technique, the pressure contours exhibit noticeable errors and zigzag patterns, indicative of the challenges faced in accurately capturing the flow behavior. In stark contrast, the pressure contours obtained with our proposed Gaussian-based method showcase a remarkable improvement, with significantly reduced errors and a smoother representation of the flow field. By effectively reducing the errors and mitigating the zigzag artifacts, our method empowers the simulations to yield more reliable and consistent results. The application of our innovative Gaussian filter-based approach has thus proven to be a transformative enhancement in accurately simulating the lid-driven cavity flow at higher Reynolds numbers. This advancement opens up new possibilities for tackling complex fluid flow problems with improved computational efficiency and precision.

\subsection{Problem-3}
\begin{figure}[htbp]
 \centering
 \vspace*{5pt}%
 \hspace*{\fill}%
\begin{subfigure}{0.50\textwidth}     
    \centering
    \includegraphics[width=\textwidth]{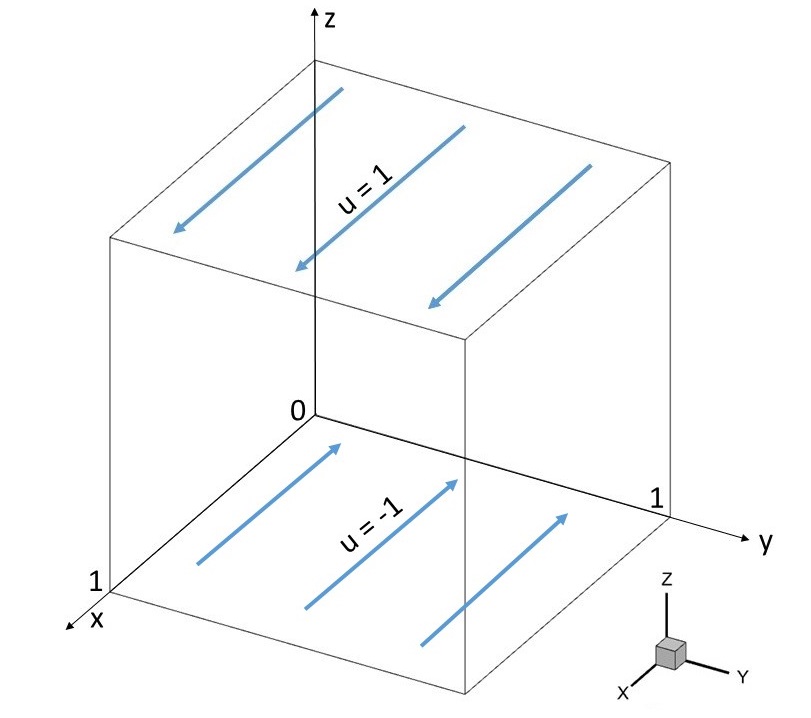}%
    \captionsetup{skip=5pt}%
    \caption{(a)}
    \label{fig:2_lid_driven_Cacity_3d}
  \end{subfigure}%
 \begin{subfigure}{0.50\textwidth}        
   \centering
    \includegraphics[width=\textwidth]{3D_GRIDS_91_91.png}%
    \captionsetup{skip=5pt}%
    \caption{(b)}
    \label{fig:3D_grids_91_91_91_c}
  \end{subfigure}
  \hspace*{\fill}

  \vspace*{8pt}%
  \hspace*{\fill}%
  \caption{(a) Illustration of the configuration in the 3D double lid-driven cavity scenario and (b) View of $91\times91\times91$ grids}
  \label{fig:2-Lid_cavity_grid}
\end{figure}
\begin{figure}
    \centering
    \includegraphics[width=\textwidth]{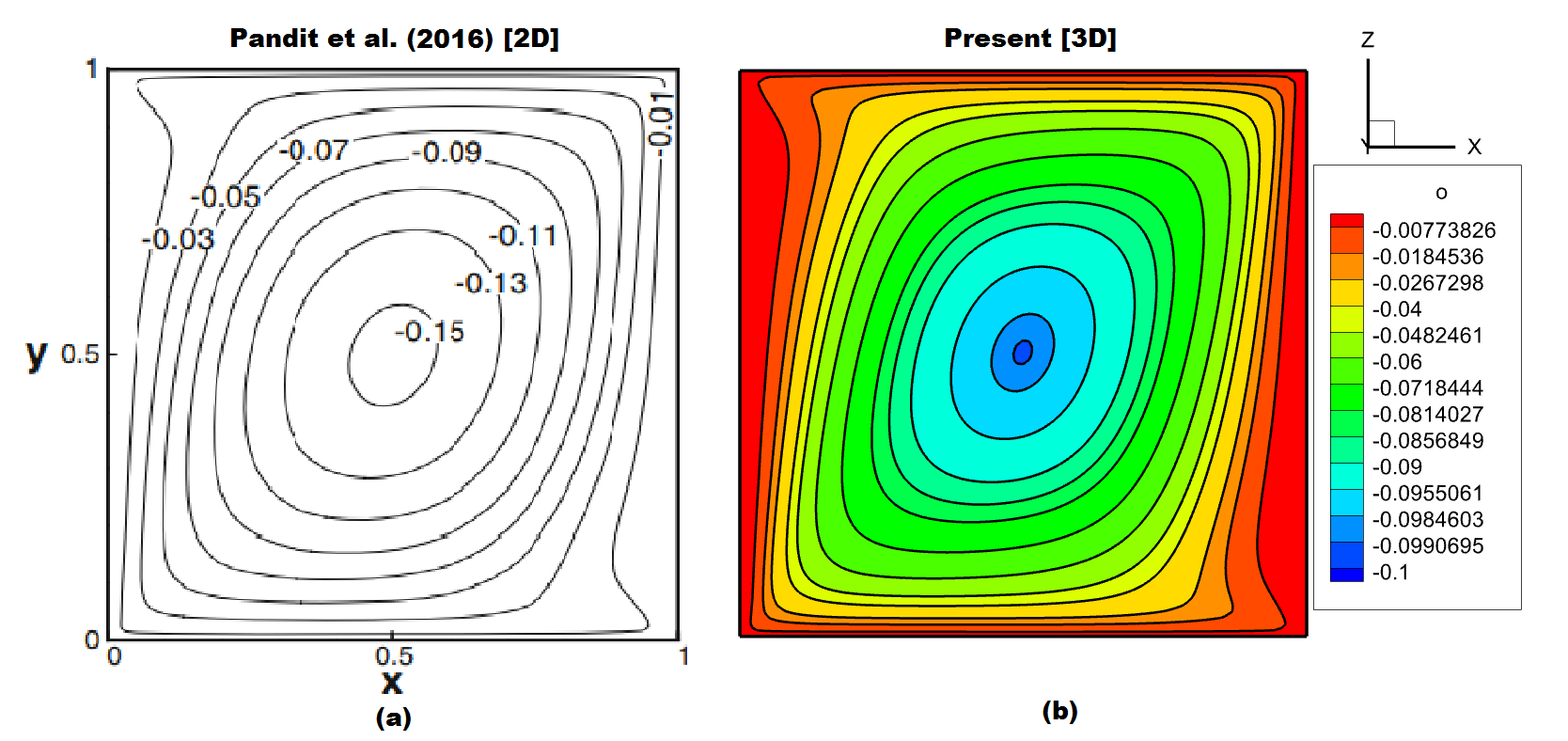}
    \caption{Comparison of the streamlines at $Re=400$ (a) Pandit et al. (2016) [2D] (b) Present [3D] at $y=0.5$ plane }
    \label{fig:double_lid_comparison}
\end{figure}

\begin{figure}[htbp]
 \centering
 \vspace*{5pt}%
 \hspace*{\fill}%
\begin{subfigure}{0.50\textwidth}     
    \centering
    \includegraphics[width=\textwidth]{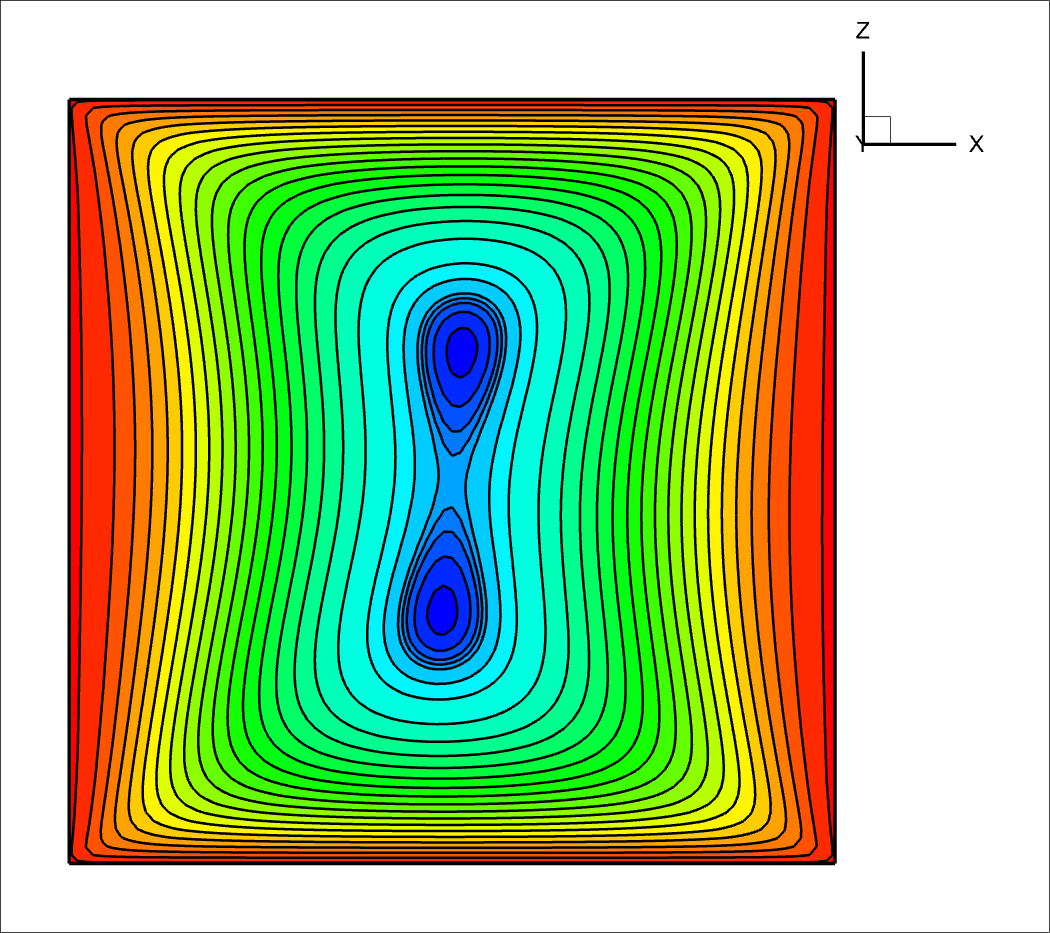}%
    \captionsetup{skip=5pt}%
    \caption{(a)}
    \label{fig:Streamlines_Re_10_double_lid}
  \end{subfigure}%
 \begin{subfigure}{0.50\textwidth}        
   \centering
    \includegraphics[width=\textwidth]{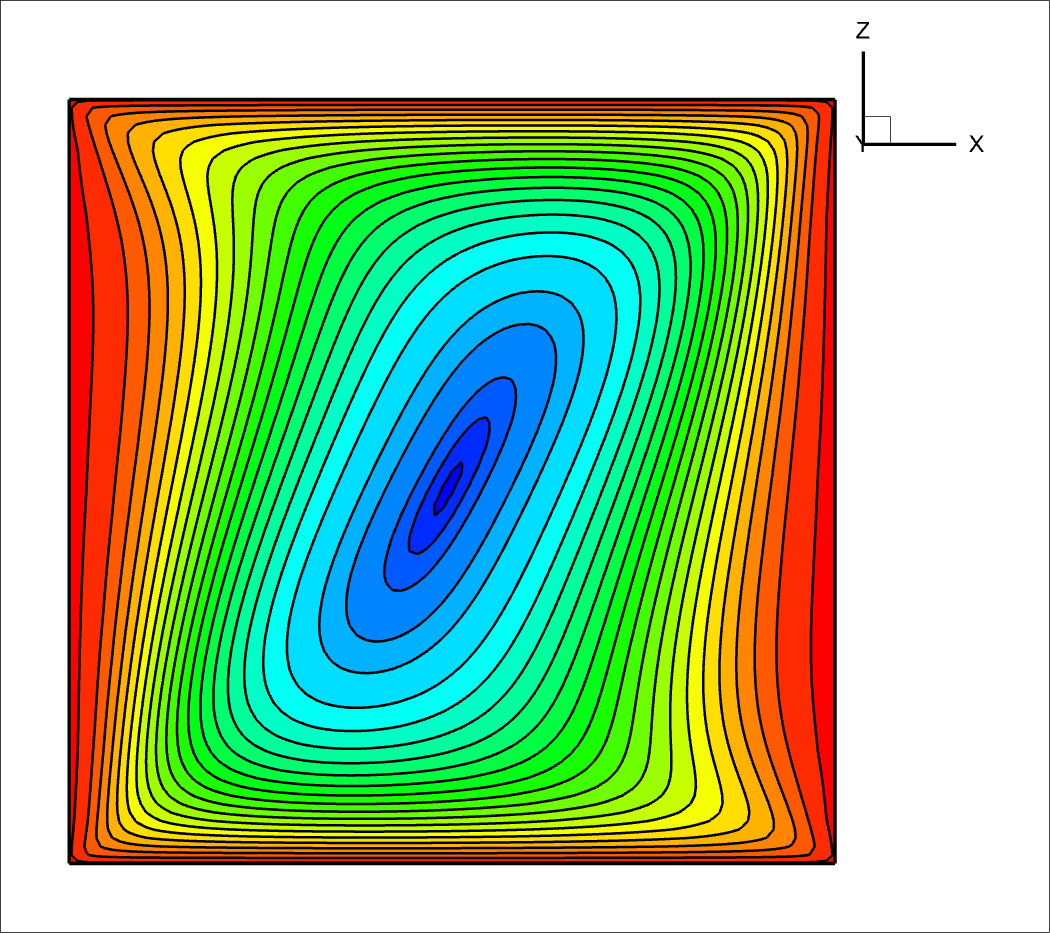}%
    \captionsetup{skip=5pt}%
    \caption{(b)}
    \label{fig:Streamlines_Re_100_double_lid}
  \end{subfigure}
  \hspace*{\fill}

  \vspace*{8pt}%
  \hspace*{\fill}%
  \begin{subfigure}{0.50\textwidth}     
    \centering
    \includegraphics[width=\textwidth]{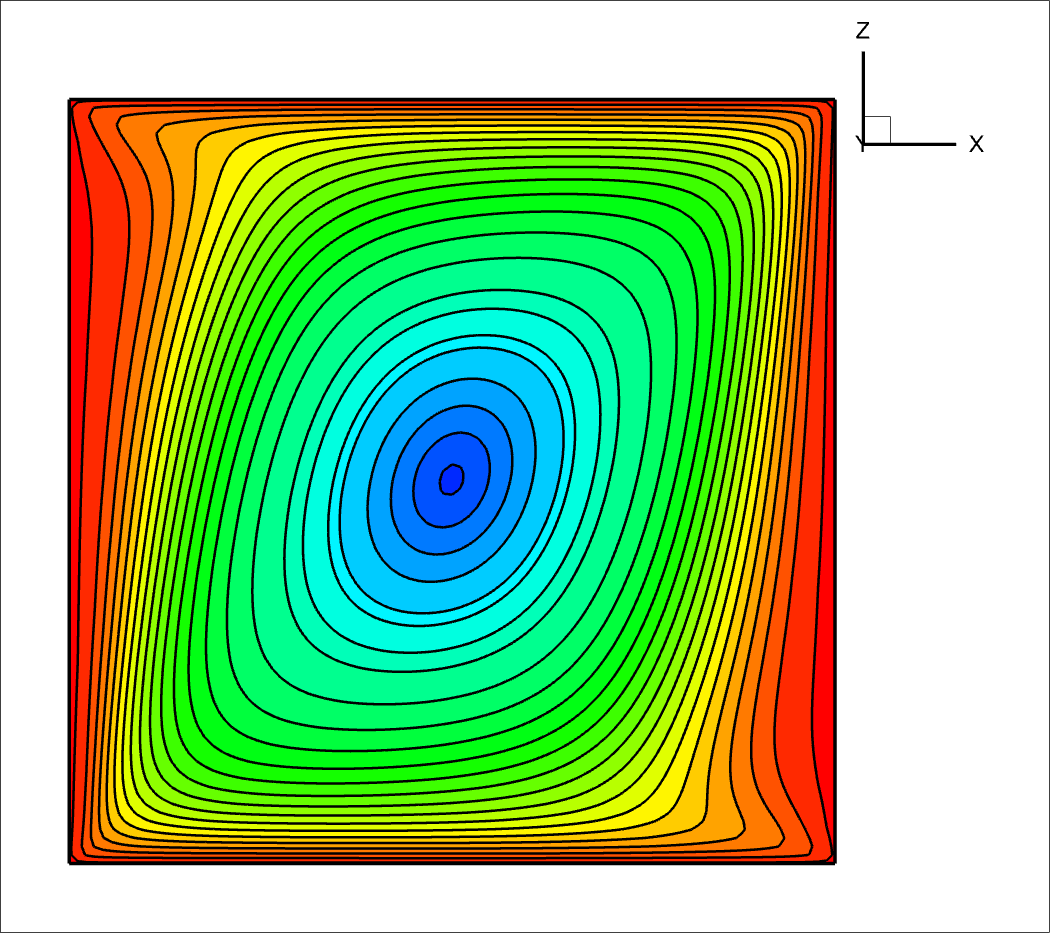}%
    \captionsetup{skip=5pt}%
    \caption{(c)}
    \label{fig:Streamlines_Re_400_double_lid}
  \end{subfigure}%
  \begin{subfigure}{0.50\textwidth}     
    \centering
    \includegraphics[width=\textwidth]{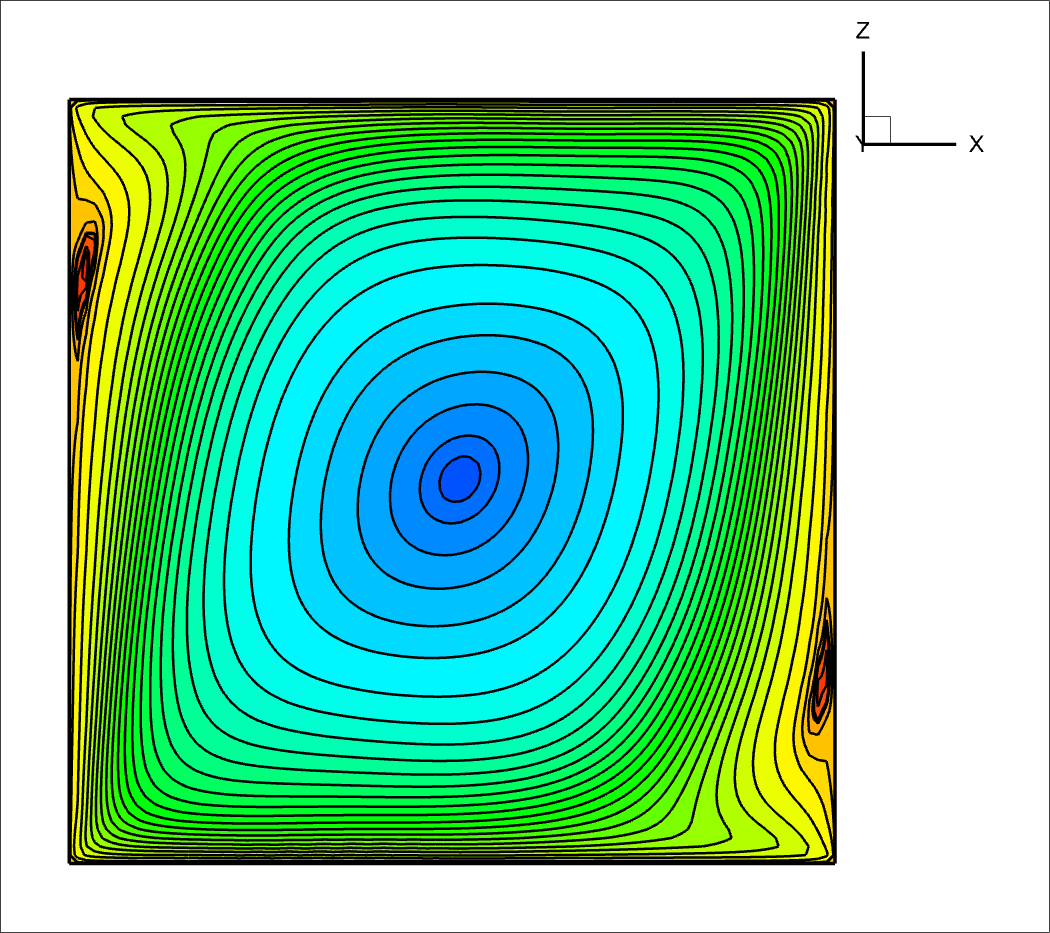}%
    \captionsetup{skip=5pt}%
    \caption{(d)}
    \label{fig:Streamlines_Re_1000_double_lid}
  \end{subfigure}%
  \caption{Streamline pattern on mid of y-plane at (a) $Re = 10$ (b) $Re = 100$ (c) $Re = 400$ and (d) $Re = 1000$ on a $91\times91\times91$ grid.}
  \label{fig:Streamlines_Double_lid_driven}
\end{figure}
\begin{figure}[htbp]
 \centering
 \vspace*{5pt}%
 \hspace*{\fill}%
\begin{subfigure}{0.50\textwidth}     
    \centering
    \includegraphics[width=\textwidth]{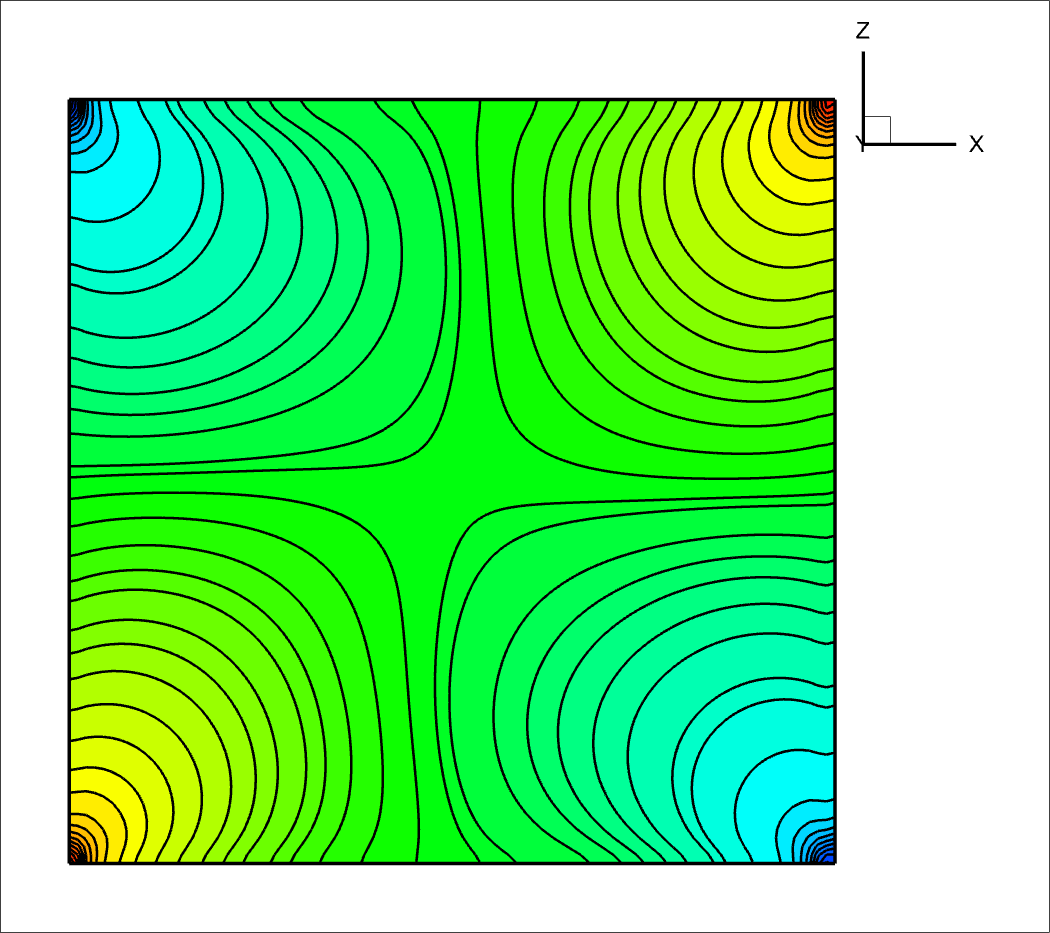}%
    \captionsetup{skip=5pt}%
    \caption{(a)}
    \label{fig:Pressure_Re_10_double_lid}
  \end{subfigure}%
 \begin{subfigure}{0.50\textwidth}        
   \centering
    \includegraphics[width=\textwidth]{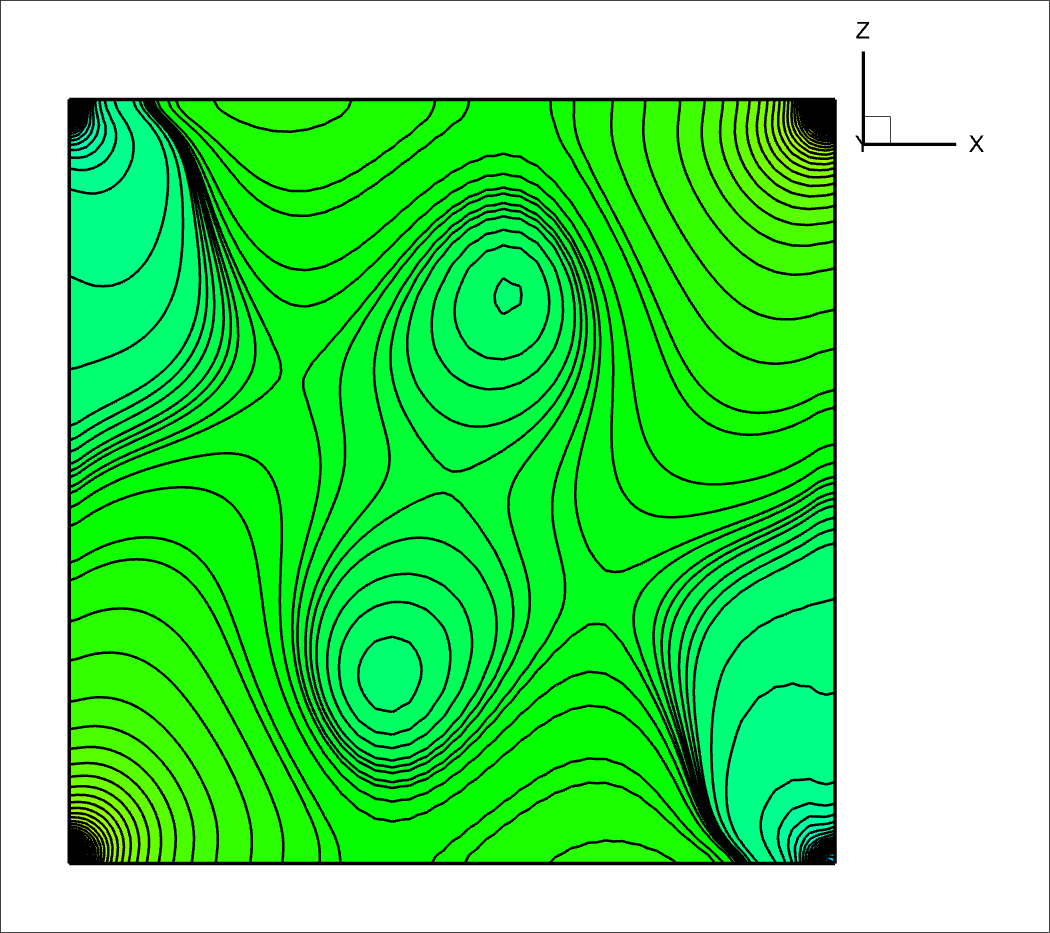}%
    \captionsetup{skip=5pt}%
    \caption{(b)}
    \label{fig:Pressure_Re_100_double_lid}
  \end{subfigure}
  \hspace*{\fill}

  \vspace*{8pt}%
  \hspace*{\fill}%
  \begin{subfigure}{0.50\textwidth}     
    \centering
    \includegraphics[width=\textwidth]{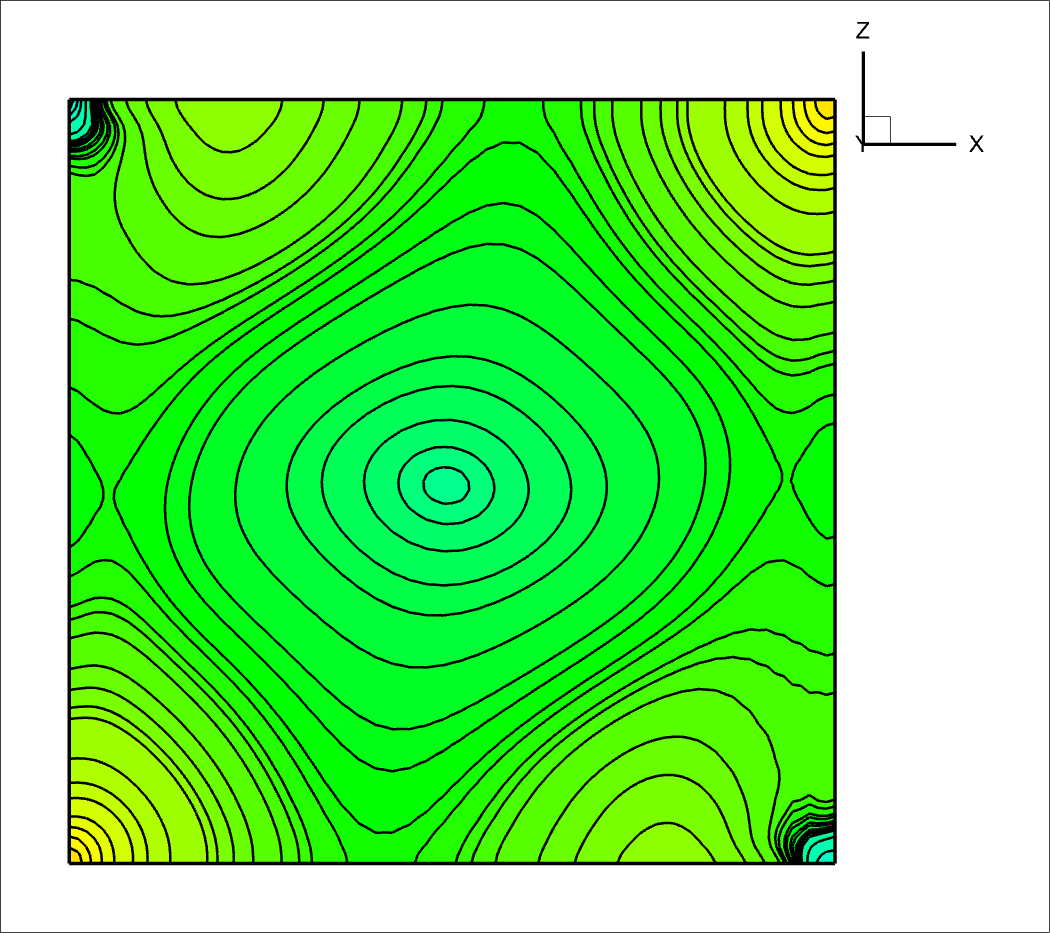}%
    \captionsetup{skip=5pt}%
    \caption{(c)}
    \label{fig:Pressure_Re_400_double_lid}
  \end{subfigure}%
  \begin{subfigure}{0.50\textwidth}     
    \centering
    \includegraphics[width=\textwidth]{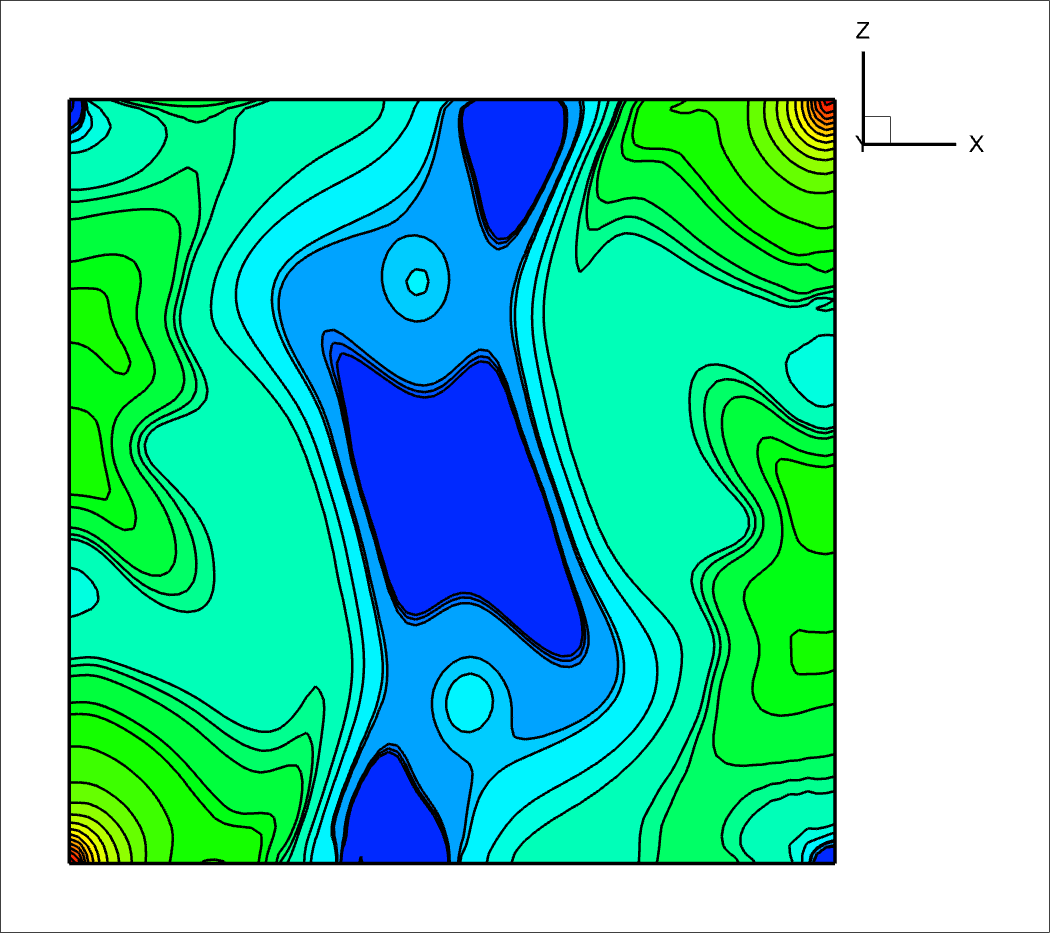}%
    \captionsetup{skip=5pt}%
    \caption{(d)}
    \label{fig:Pressure_Re_1000_double_lid}
  \end{subfigure}%
  \caption{Pressure distribution on the middle of y-plane ($y=0.5$) at (a) $Re = 10$, (b) $Re = 100$, (c) $Re = 400$, and (d) $Re = 1000$ on $91\times91\times91$ grid.}
  \label{fig:Pressure_Double_lid_driven}
\end{figure}

In the final problem, we investigate the dynamics within a cubical cavity driven by dual lids using a super compact higher-order scheme that incorporates our proposed pressure correction technique. While this problem has been extensively studied in its 2D \cite{Blohm_2002, Perumal_2010, Pandit_2016} form, there exists a scarcity of literature that addresses its 3D counterpart. By implementing the super compact scheme and integrating our novel pressure handling technique, we aim to enhance the accuracy and efficiency of the numerical simulations. The super compact scheme is renowned for its high order of accuracy and superior performance in capturing finer details of the flow dynamics. By combining this advanced numerical scheme with our innovative pressure-correction technique, we seek to address the challenges posed by the double lid-driven cavity flow problem. Our focus rests on enhancing both accuracy and efficiency. Within this investigation, we performed simulations across various Reynolds numbers: specifically, $Re = 10, 100, 400$, and $1000$. Figure \ref{fig:2-Lid_cavity_grid} presents the schematic diagram of the problem. As illustrated in the diagram, the two parallel lids are set in motion, but in opposing directions. The remaining boundary and initial conditions are kept consistent with those used in the previous single lid-driven cavity problem. Across the upper face of the cubic cavity, a uniform velocity is specified along the $x$ direction. Conversely, the same negative uniform velocity is applied to the lower face of the cavity, as shown in Figure \ref{fig:2-Lid_cavity_grid}.\\
In Figure \ref{fig:double_lid_comparison}, a comparison is depicted between the streamlines of the 2D lid-driven cavity\cite{Pandit_2016} and the current 3D simulation results. This comparison is conducted at a fixed plane $y=0.5$ with a Reynolds number, $Re=400$. The visual patterns of the streamlines exhibit close resemblance, although minor quantitative differences are observed due to the influence of the three-dimensional effects. Figure \ref{fig:Streamlines_Double_lid_driven} illustrates the streamlines pattern at the middle of the $y$-plane for all Reynolds numbers considered. Indeed, the steady-state flow behavior in the two parallel lid-driven cubical cavity exhibits interesting variations with changing Reynolds numbers. We noted the presence of two primary vortex centers within the flow pattern for the case of $Re=10$. As Reynolds number increases, fluid flow becomes more structured and the primary vortex center becomes more distinct and predominant. Specifically, for $Re$ 100, and 400, the flow is characterized by the formation of a single primary vortex within the cavity. This primary vortex dominates the flow pattern, and no other significant vortices are observed.\\
However, as we increase the Reynolds number to $Re = 1000$, the flow behavior undergoes a notable transformation. At this higher Reynolds number, the flow becomes more dynamic, leading to the appearance of additional secondary vortices. In particular, two corner vortices become prominent in the left upper and bottom right corners of the cavity. These secondary vortices are indicative of the enhanced complexity and increased circulation within the flow field. The presence of two corner vortices at $Re = 1000$ demonstrates the influence of higher Reynolds numbers on the flow characteristics and highlights the non-linear nature of the flow behavior. It also indicates the significant impact of the proposed numerical scheme and pressure handling technique in accurately capturing such intricate flow phenomena, which is crucial for understanding and predicting fluid dynamics in practical applications.\\
Figure \ref{fig:Pressure_Double_lid_driven} provides valuable insights into the pressure distribution within the two-sided lid-driven cubical cavity for different Reynolds numbers. The pressure contours are visualized in the middle of the y-plane, allowing us to observe the spatial variations in pressure and understand the impact of Reynolds numbers on the flow pattern. At lower Reynolds numbers ($Re$ = 10, 100 and 400), the pressure contours exhibit a relatively smooth and uniform distribution across the cavity. The pressure values are relatively low, indicating a more laminar and less turbulent flow behavior. The pressure contours form coherent patterns that reflect the dominant primary vortex present in the flow at these Reynolds numbers. However, as we move to a higher Reynolds number ($Re =  1000$), the pressure contours exhibit significant changes. The most notable effect is the emergence of irregularities and higher pressure gradients in the flow field. These irregularities are the result of increased turbulence and the presence of additional secondary vortices at $Re = 1000$, as discussed earlier. The two corner vortices at $Re = 1000$ introduce complexities in the flow as displayed in the Figure \ref{fig:Pressure_Double_lid_driven}.

\begin{figure}[htbp]
 \centering
 \vspace*{5pt}%
 \hspace*{\fill}%
\begin{subfigure}{0.50\textwidth}     
    \centering
    \includegraphics[width=\textwidth]{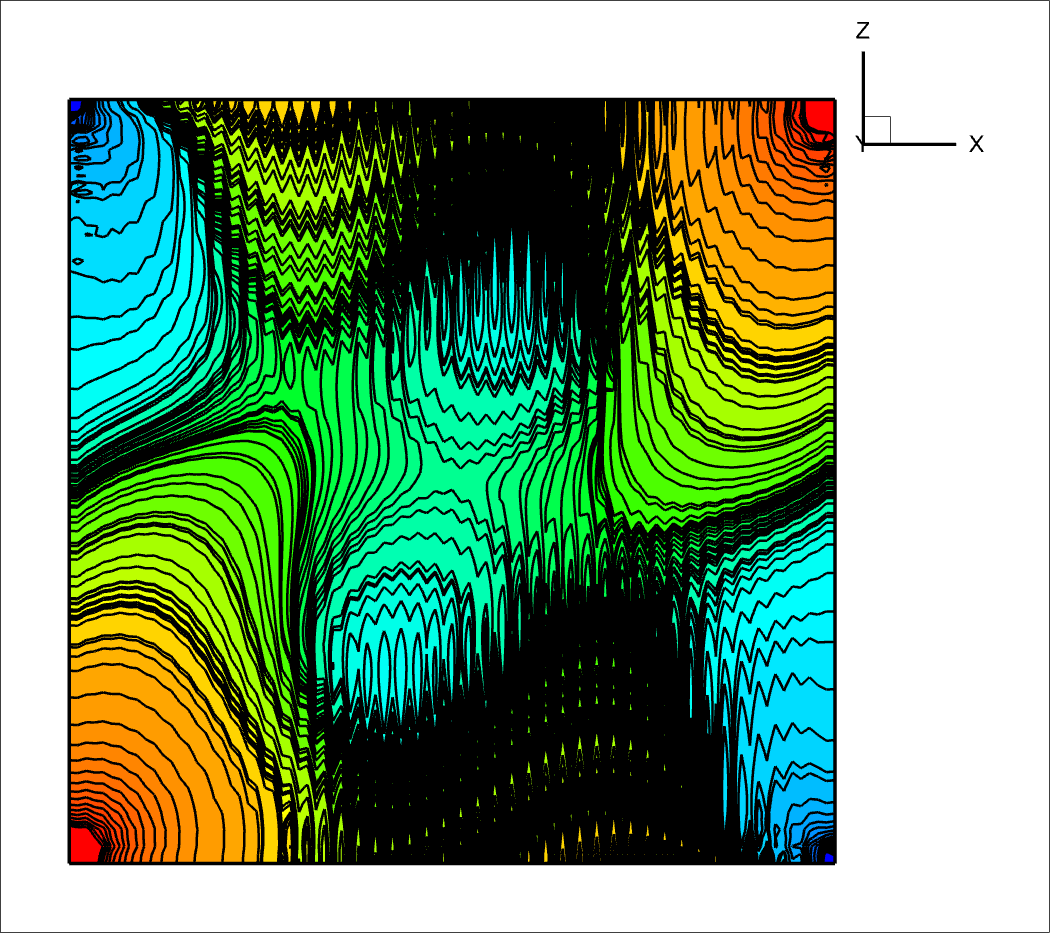}%
    \captionsetup{skip=5pt}%
    \caption{(a)}
    \label{fig:Pressure_Re_100_zig_zag_double_lid}
  \end{subfigure}%
 \begin{subfigure}{0.50\textwidth}        
   \centering
    \includegraphics[width=\textwidth]{Pressure_Re_100_double_lid.png}%
    \captionsetup{skip=5pt}%
    \caption{(b)}
    \label{fig:Pressure_Re_100_double_lidc}
  \end{subfigure}
  \hspace*{\fill}

  \vspace*{8pt}%
  \hspace*{\fill}%
  \begin{subfigure}{0.50\textwidth}     
    \centering
    \includegraphics[width=\textwidth]{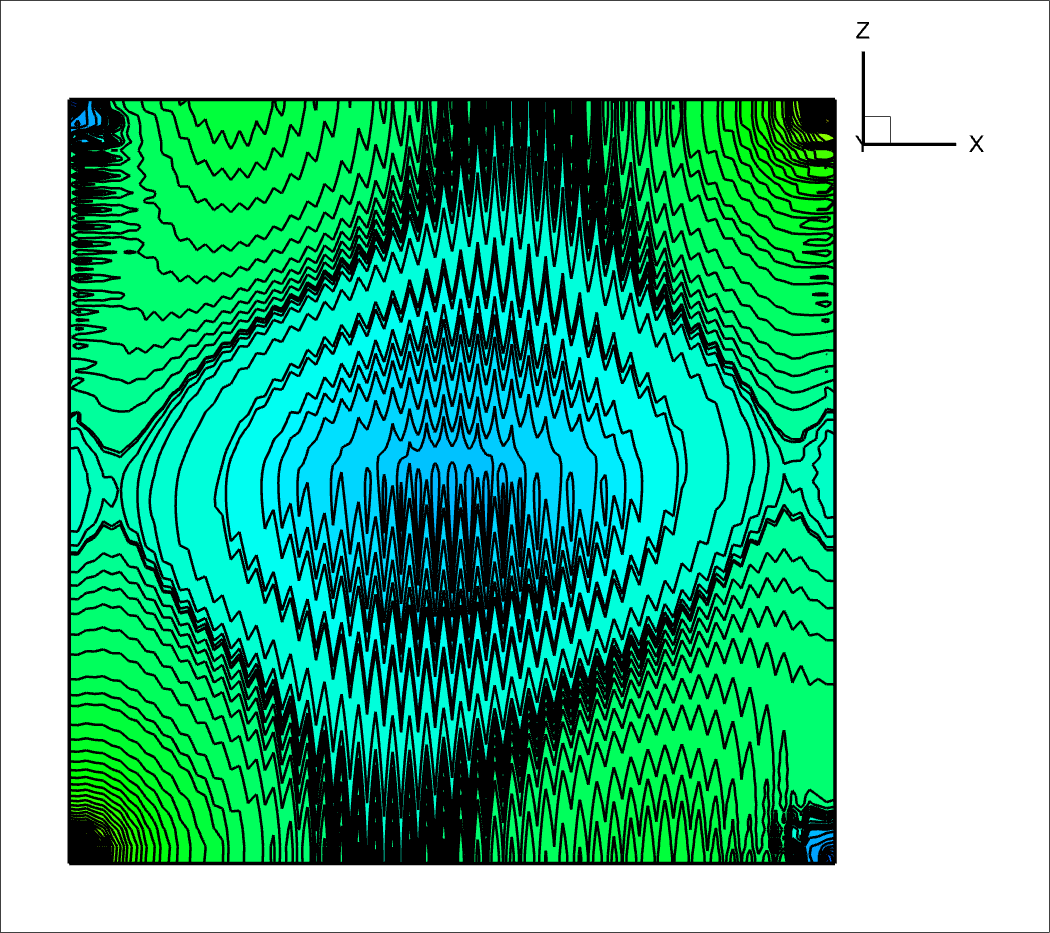}%
    \captionsetup{skip=5pt}%
    \caption{(c)}
    \label{fig:Pressure_Re_400_Zig_zag_double_lid}
  \end{subfigure}%
  \begin{subfigure}{0.50\textwidth}     
    \centering
    \includegraphics[width=\textwidth]{Pressure_Re_400_double_lid.png}%
    \captionsetup{skip=5pt}%
    \caption{(d)}
    \label{fig:Pressure_Re_400_double_lidc}
  \end{subfigure}%
  \caption{Comparison of the pressure contours (a,c) without using (left) and (b,d) with using (right) our proposed pressure correction technique at the same number of pressure iterations for (a,b) $Re=100$, (c,d) $Re=400$.}
  \label{fig:Pressure_zig_zag_double}
\end{figure}

{\small\begin{table}[htbp]
\caption{\small Comparative Analysis of Various Parameters in the 3D Double Lid-Driven Cavity Problem-3}\label{Pressure_Iterations_2}
\centering
 \begin{tabular}{cccccc}  \hline \hline
  &    &     Without Using Proposed Pressure Correction Technique  &  &  &       \\ \hline
& $Re$    &   $\lambda$    &  $\Delta$t  &  Iterations  &    \\ \hline 
&10          &  0.1  &    0.01  &    15\\
&100         &  0.06   & 0.005  &    24\\
&400             &  0.02  &  0.001 &    37    \\
&1000             &  0.02  &  0.0001 &    54    \\\hline 
&    &     Using Proposed Pressure Correction Technique &    &    &    \\ \hline 
&10          &  0.1  &    0.02  &    3\\
&100         &  0.06   & 0.005  &    4\\
&400             &  0.02  &  0.002 &    4    \\
&1000             &  0.02  &  0.001 &    5   \\
\hline\hline
 \end{tabular}
\end{table}
}
We have generated the results both with and without the implementation of the proposed pressure correction technique. Table \ref{Pressure_Iterations_2} presents the comparison of the computational parameters used for solving the double lid-driven cavity problem, both with and without the application of our proposed technique. The table presents the $\Delta t$ values and the corresponding number of iterations needed for pressure calculations at various Reynolds numbers. Pressure iterations are terminated when $|\nabla \cdot \mathbf{v}|_{\max}$ reaches a tolerance limit of $10^{-3}$. The results are obtained utilizing a grid size of $91 \times 91 \times 91$. Interestingly, the results highlight a notable decrease in the required iterations with the adoption of our proposed pressure correction technique. The decrease in iterations directly corresponds to a significant reduction in computational cost, rendering the proposed technique remarkably efficient in solving this problem too. The conventional approach necessitates a time step size of $\Delta t = 0.0001$ for the convergence of the iterative method, whereas with the utilization of our pressure correction technique, the required time step is only $\Delta t = 0.001$. At the same time, our technique takes only 5 inner iterations to converge to the desired accuracy in each time step, whereas the conventional approach needs 54 inner iterations for the same task. This also demonstrates the efficiency and computational cost-effectiveness of the proposed pressure correction technique. Once more, this demonstrates the efficacy of the proposed pressure correction method.

\begin{figure}[htbp]
 \centering
 \vspace*{5pt}%
 \hspace*{\fill}%
\begin{subfigure}{0.50\textwidth}     
    \centering
    \includegraphics[width=\textwidth]{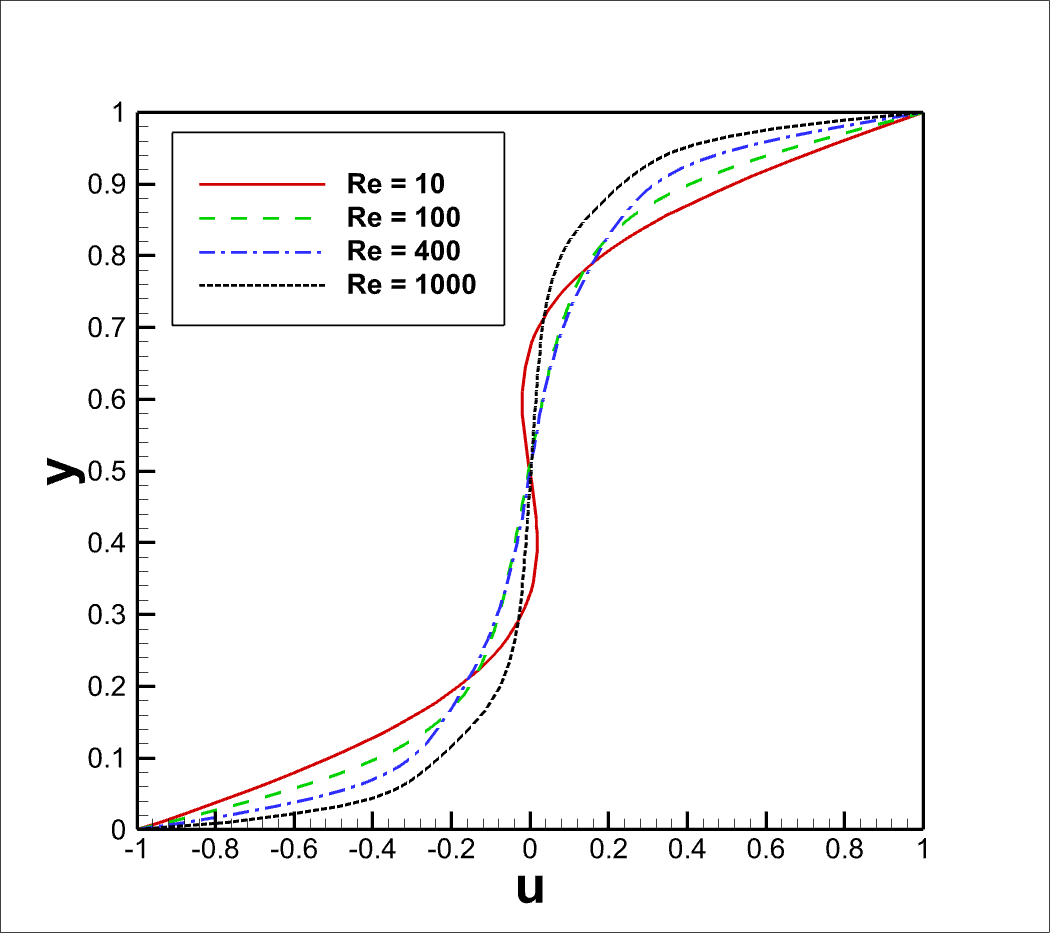}%
    \captionsetup{skip=5pt}%
    \caption{(a)}
    \label{fig:u_y_double_lid_driven}
  \end{subfigure}%
 \begin{subfigure}{0.50\textwidth}        
   \centering
    \includegraphics[width=\textwidth]{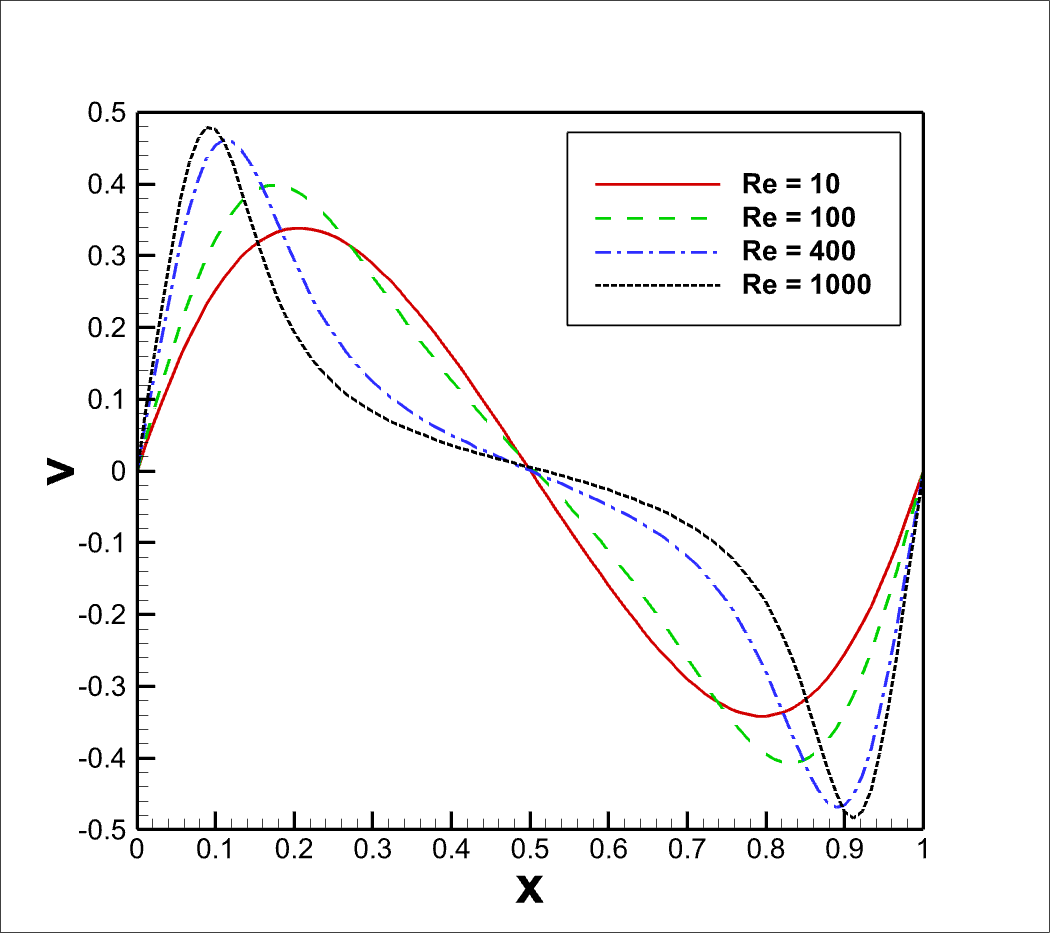}%
    \captionsetup{skip=5pt}%
    \caption{(b)}
    \label{fig:x_v_double_lid_driven}
  \end{subfigure}
  \hspace*{\fill}

  \vspace*{8pt}%
  \hspace*{\fill}%
  \caption{ Centerline velocity variation for different Reynold number (a) Horizontal velocity profile along the vertical centerline, (b) Vertical velocity profile along the horizontal centerline}
  \label{fig:velocity_variation}
\end{figure}

In Figure \ref{fig:Pressure_zig_zag_double}, we present a comparison of the pressure distributions for $Re=100, 400$ obtained with and without the utilization of our proposed pressure correction technique, using the same number of pressure iterations. The striking contrast between the two sets of results serves as a clear testament to the effectiveness of our proposed scheme. When we examine the pressure contours obtained without implementing our technique, we can observe several undesirable artifacts in the form of errors and zigzag patterns. These imperfections indicate that the traditional pressure handling method is not as effective for low-pressure iterations in capturing the true pressure distribution and resolving the flow characteristics accurately.\\
In Figure \ref{fig:velocity_variation}, we illustrate a comprehensive comparison of the horizontal velocity profiles along the vertical centerline, and the vertical velocity profile along the horizontal centerline within the square cavity. The simulations are conducted for various Reynolds numbers, i.e., $Re = 10, 100, 400,$ and $1000$. The objective of this analysis is to understand the impact of increasing Reynolds numbers on the flow dynamics and velocity profiles within the cavity. By examining the variations in horizontal and vertical velocities with increasing Reynolds numbers, we can gain valuable insights into the flow. 
Specifically, we can observe that the peak $v$-velocity values increase with the increasing Reynolds number. This behavior is expected as higher Reynolds numbers lead to more energetic and complex flow patterns, resulting in greater vortical motion and increased $v$-velocity magnitudes. However, the behavior of the $u$-velocity is more intriguing. As the Reynolds number increases, we notice a distinct difference in the $u$-velocity profiles depending on the location along the y-axis. When $y < 0.5$ (i.e., the bottom half of the cavity), the $u$-velocity increases with the increasing  Reynolds number. This phenomenon can be attributed to the influence of the imposed lid-driven motion, which propels the fluid in the $x$-direction more vigorously as the Reynolds number increases. Consequently, the $u$-velocity magnitudes intensify in the bottom half of the cavity. On the other hand, when $y > 0.5$ (i.e., the top half of the cavity), we observe a different trend in the u-velocity profiles. Here, the $u$-velocity values decrease with the increasing Reynolds number. This behavior arises due to the interplay between the imposed lid-driven flow and the developing vortical structures in the cavity. As the Reynolds number increases, the vortices become more pronounced and influence the fluid motion in the top half of the cavity, resulting in a reduction of the u-velocity magnitudes in this region.

\section{Conclusion}
In conclusion, we can say that we have introduced a novel and effective technique for pressure correction, which has demonstrated its ability to significantly reduce the computational cost of pressure calculation, and overall numerical simulation. Our approach is applied to three different problems: the 3D Burger's equation and two variations of the lid-driven cavity problem.
In the case of the 3D Burger's equation, we have successfully obtained fourth-order accurate solutions, and our results have shown good agreement with the analytical solutions. For the lid-driven cavity problems (both single and double lid-driven cases), which serve as a widely-used benchmark for assessing numerical methods, we have conducted extensive simulations at various Reynolds numbers ($Re$=10, 100, 400, 1000). Our computed results are compared to existing literature, and we have observed excellent agreement. The simulations have been conducted using different Reynold numbers, and our proposed technique has consistently shown superior performance in reducing computational costs while maintaining accuracy. Furthermore, we have investigated the flow patterns and pressure contours for different Reynolds numbers, revealing the distinct effects of fluid behavior as $Re$ increases. Our technique has successfully captured these flow phenomena and demonstrated its capability to handle complex flow patterns effectively. In summary, our proposed technique has proven to be a valuable addition to the field of computational fluid dynamics. By significantly reducing computational costs without compromising accuracy, it opens up new possibilities for conducting simulations of more extensive and complex fluid flow problems. The computational efficiency gained through our method allows for faster convergence rates and enables the exploration of a broader range of Reynolds numbers and grid resolutions. As we move forward, the versatility and robustness of our proposed technique make it applicable to diverse scenarios in fluid dynamics research and engineering simulations. Its simplicity and effectiveness pave the way for its adoption in a wide range of practical scenarios.\\
Overall, this study demonstrates the potential of our new pressure correction technique in advancing the field of numerical simulations for three-dimensional fluid flow problems. The successful application to both analytical and benchmark problems showcases the reliability of our approach and encourages further exploration and development in this direction.\vspace{20pt}\\
{\large\textbf{Data Availability}}\vspace{5pt}\\
The data that support the findings of this study are available from the corresponding author upon reasonable request.\vspace{15pt}\\
{\large\textbf{Compliance with Ethical Standards}}\vspace{5pt}\\
\textbf{Conflict of Interest:} All authors declare that they have no conflict of interest.\vspace{5pt}\\
\textbf{Ethical approval:} This article does not contain any studies with human participants or animals performed by any of the authors.\\

\end{document}